\documentclass[11pt,leqno]{article}
\usepackage{graphicx}
\usepackage{amsfonts}
\usepackage{amsmath}
\usepackage{amsthm}
\usepackage{rotating}
\usepackage{verbatim}
\usepackage[export]{adjustbox}
\usepackage{varioref}
\usepackage{float}

\usepackage{color}
\usepackage{subcaption}

\usepackage{ amssymb }
\usepackage{hyperref}
\usepackage[round]{natbib}

\newcommand\balpha{\mbox{\boldmath${\alpha}$}}

\newcommand\vep{{\varepsilon}}

\newcommand\bA{{\bf A}}

\newcommand\mbK{{\mathbb K}}
\newcommand\mcL{{\mathcal L}}

\newcommand\bY{{\bf Y}}

\usepackage{enumerate}

\newcommand\mbR{{\mathbb R}}

\DeclareMathOperator\E{E}
\DeclareMathOperator\Span{span}

\newtheorem{lemma}{Lemma}
\newtheorem{theorem}{Theorem}

\theoremstyle{definition}

\newtheorem{asu}{Assumption}
\newcounter{subassumption}[asu]
\renewcommand{\thesubassumption}{(\textit{\roman{subassumption}})}
\makeatletter
\renewcommand{\p@subassumption}{\theasu}
\makeatother
\newcommand{\subasu}{
  \refstepcounter{subassumption}%
  \thesubassumption~\ignorespaces}

\begin{document}
\title{Optimal Prediction for Additive Function-on-Function Regression}

\author{Matthew Reimherr\footnote{Corresponding author: Matthew Reimherr, 411 Thomas Building, University Park, PA 16802, \href{mreimherr@psu.edu}{mreimherr@psu.edu}} \ and Bharath Sriperumbudur \\ Pennsylvania State University 
\and Bahaeddine Taoufik \\
Lynchburg College 
}

\date{}
\maketitle

\begin{abstract}
As with classic statistics, functional regression models are invaluable in the analysis of functional data.  While there are now extensive tools with accompanying theory available for linear models, there is still a great deal of work to be done concerning nonlinear models for functional data.  In this work we consider the Additive Function-on-Function Regression model, a type of nonlinear model that uses an additive relationship between the functional outcome and functional covariate.  We present an estimation methodology built upon Reproducing Kernel Hilbert Spaces, and establish optimal rates of convergence for our estimates in terms of prediction error.  We also discuss computational challenges that arise with such complex models, developing a representer theorem for our estimate as well as a more practical and computationally efficient approximation.  Simulations and an application to Cumulative Intraday Returns around the 2008 financial crisis are also  provided.
\end{abstract}

\section{Introduction}
Functional data analysis (FDA) concerns the statistical analysis of data where 
one of the variables of interest is a function. FDA has seen rapidly increasing 
interest over the last few decades and has successfully been applied to a 
variety of fields, including economics, finance, the geosciences, and the health 
sciences. One of the most fundamental tools in statistics is linear 
regression, as such, it has been a major area of research in FDA. While the 
literature is too vast to cover here, we refer readers to \citet{ramsay2006functional,ramsay2009functional,horvath2012inference,KRBook}, which provide introductions to FDA, as well as \citet{morris2015functional}, which provides a broad overview of methods for functional linear regression.

A major challenge of functional regression is handling functional predictors. 
At least conceptually, a functional predictor means having a large number 
(theoretically infinite) of predictors that are all highly correlated. To 
handle such a setting, certain regularity conditions are imposed to make the 
problem tractable. Most of these conditions are directly or indirectly related 
to the smoothness of the parameter being estimated. However, the convergence 
rates of the resulting estimators then depend heavily on these assumptions, and 
the rates are not parametric when the predictor is infinite dimensional.

One of the most well studied models in FDA is the functional linear model. 
Commonly, one distinguishes between function-on-scalar, scalar-on-function, and 
function-on-function regression when discussing such models, with first term denoting the type of response and the second term denoting the type of covariate. The convergence 
rates for function-on-scalar regression are usually much faster than for the 
scalar-on-function or function-on-function. Methodological, theoretical, and 
computational issues related to functional linear models are now well 
understood. More recently, there has been a growing interest in developing 
nonlinear regression models. While it is natural to begin examining nonlinear 
models after establishing the framework for linear ones, there is also a 
practical need for such models. 
Functional data may contain complicated temporal dynamics, which may exhibit 
nonlinear patterns that are not well modeled assuming linearity; 
\citet{fan2015functional} examine this issue deeply. 

Nonlinear 
regression methods for FDA have received a fair amount of attention for the 
scalar-on-function setting, while function-on-function regression models, where 
the relationship between the response and covariates is believed to be 
nonlinear, have received considerably less attention. Concerning nonlinear 
scalar-on-function regression, \citet{james2005functional} introduced a 
functional single index model, where the outcome is related to a linear 
functional of the predictor through a nonlinear transformation. This work would 
later be extended in \citet{fan2015functional}, allowing for a potentially 
high-dimensional number of a functional predictors. \citet{preda2007regression} explored fitting a fully nonlinear model using 
reproducing kernel Hilbert spaces (RKHS). In contrast, \citet{muller2013continuously} simplified the form of the 
nonlinear relationship by introducing the functional additive model, which 
combines ideas from functional linear models and scalar additive models 
\citep{hastie1990generalized}. Optimal convergence rates for the functional additive model were then 
established by \citet{wang2015optimal}, which generalized the work of 
\citet{cai2012minimax} in the linear case. An alternative to the functional 
additive model was given in \citet{zhu2014structured} who first expressed the 
functional predictor using functional principal components analysis, FPCA, 
and then built an additive model between the 
outcome and scores.  An extension to generalized linear models can be found in \cite{mclean2014,du2014penalized}.   

Moving to function-on-function regression, \citet{lian2007nonlinear} extended 
the work of \citet{preda2007regression} to functional outcomes, which was then also considered in 
\citet{kadri2010nonlinear}.
Most relevant to the present paper is the work of \citet{scheipl2015functional}
who extended the work of 
\citet{muller2013continuously} by introducing an additive model for function-on-function regression.  They used a general trivariate tensor product basis approach for estimation, which allowed them to rely on  \textit{GAM} from the \textit{MGCV} package in {\tt R} to carry out the computation, as is implemented in the \textit{Refund} package.  \citet{ma2016continuously}, examining the same model, considered a binning estimation technique combined with FPCA.  In addition, they were able to prove convergence of their estimators, but made no mention of optimality while also needing a great deal of assumptions which are challenging to interpret.  Another estimation technique was examined in \citet{kim2017additive}, which was similar to the trivariate tensor product approach of \citet{scheipl2015functional}, but two of the bases are explicitly assumed to be orthogonal B-splines, while the third comes from an FPCA expansion.  However, as with \citet{scheipl2015functional}, no theoretical justification is provided.  Lastly, in very recent work, \cite{sun2017optimal} considered the case of using an RKHS framework to estimate a function-on-function linear model.  Extending the work the \cite{cai2012minimax}, they were able to establish the optimality of their procedure.  Our work can be viewed as extending this work to nonlinear relationships via a function-on-function additive model.

The goal of this work is to develop a penalized regression framework based on 
Reproducing Kernel Hilbert Spaces, RKHS, for fitting the additive 
function-on-function regression model, AFFR \citep{scheipl2015functional}. A major 
contribution of this work is to provide optimal convergence rates of our 
estimators in terms of prediction error, and that this rate is the same as for the scalar outcome setting \citep{wang2015optimal}.  We also discuss computational aspects of our approach, as the RKHS structure allows for a fairly efficient computation as compared to the trivariate tensor product bases that have been used previously.  Background and the model are introduced in Section \ref{s:model}.  Computation is discussed in Section \ref{s:comp}, while theory is presented in Section \ref{s:theory}.  We conclude with a numeric study consisting of simulations and an application to financial data.

\section{Model and Background}\label{s:model}
We assume that we observe i.i.d pairs $\{ (X_i(t), Y_i(t)): i= 1, \dots, n, \ t 
\in [0,1]\}$. The functions could be observed on other intervals, but as long 
as they are closed and bounded, then they can always be rescaled to be $[0,1]$, 
thus it is common in FDA to work on the unit interval. Both the outcome, $Y_i(t)$, and 
$X_i(t)$ are assumed to be completely observed functions, a practice sometimes 
referred to as \textit{dense functional data analysis} \citep{KRBook}; practically this means that the curve reconstruction contributes a comparatively small amount of uncertainty to the final parameter estimates.  More rigorous definitions can be found in \citet{cai2011optimal,li2010uniform,zhang2016sparse}.  For sparsely observed curves, it is usually better to use more tailored approaches such as PACE \citep{yao2005functional}, FACE  \citep{xiao2017fast}, or MISFIT \citep{petrovich2018functional}.

The additive 
function-on-function regression model is defined as
\[
Y_i(t) = \int_0^1 g(t,s, X_i(s))\, ds + \vep_i(t).
\]
We assume that the functions $X_i$, $\vep_i$, and $Y_i$ are elements of 
$L^2[0,1]$, which is a real separable Hilbert space. The trivariate function, 
$g(t,s,x)$ is assumed to be an element of an RKHS, $\mbK$.

Recall that an RKHS is a Hilbert space that possesses the \textit{reproducing 
property}, namely, we assume that $\mbK$ is a Hilbert space of functions from 
$[0,1] \times [0,1] \times \mbR \to \mbR$, and that there exists a kernel 
function $k(t,s,x, t', s', x')= k_{t,s,x}(t',s',x')$ that satisfies
\[
f(t,s,x) = \langle k_{t,s,x}, f \rangle_{\mbK},
\]
for any $f \in \mbK$. There is a one-to-one correspondence between $\mbK$ and 
$k$, thus choosing the kernel function completely determines the resulting RKHS. 
 The functions in $\mbK$ inherit properties from $k$, in particular, one can 
choose $k$ so that the functions in $\mbK$ possess some number of derivatives, 
or satisfy some boundary conditions. In addition, many Sobolev spaces, which 
are commonly used to enforce smoothness conditions, are also RKHS's. We refer an 
interested reader to \citet{berlinet2011reproducing} for further details.


We propose to estimate $g$ by minimizing the following penalized objective:
\begin{align}
RSS_\lambda(g) = \sum_{i=1}^{n} \int_0^{ 1}\left (Y_{i}(t)- \int_0^{ 1} 
g(t,s,X_{i}(s))\, ds \right)^2 dt+ \lambda \| g \|^{2} _{\mbK}, \label{e:rss}
\end{align}
i.e., $$\hat g = \arg\inf_{g\in\mbK} RSS_\lambda(g),$$ where $\lambda>0$. As we will see in the next section, an explicit solution to this minimization problem exists due to the reproducing property.  However, we will also discuss using FPCA to help reduce the computational burden.

\section{Computation} \label{s:comp}
One of the benefits of using RKHS methods is that one can often get an exact 
solution to the corresponding minimization problem such as the one in \eqref{e:rss}, due to the representer theorem \citep{Kimeldorf-71}. This also turns out to be 
the case here, however, later on we will discuss using a slightly modified version 
that still works well and is easier to compute.  The expression we derive is quite a bit simpler than the analogs derived in \cite{cai2012minimax,wang2015optimal,sun2017optimal}; this is partly due to our use of functional principal components, which simplify the expression and also provide an avenue for reducing the computational complexity of the problem, and also due to our use of the RKHS norm penalty when fitting the model (where as others used a more general penalty term).  

%
Using the reproducing property we have 
\[
 \langle k_{t,s, X_{i}(s)}, g \rangle_{\mbK} = g(t,s,X_{i}(s)) \hspace 
{0.4cm}\text{for} \hspace {0.4cm} i=1,2,...,n.
\]
We then have that 
 \begin{equation}
 \int_0^{ 1} g(t,s,X_{i}(s)) ds = \int_0^{ 1}\langle g, k_{t,s, X_{i}(s)} 
\rangle_{\mbK} ds= \left \langle g, \int_0^{ 1} k_{t,s, X_{i}(s)} ds \right 
\rangle_{\mbK},
\end{equation}
which is justified by the integrability constraints inherent in Assumption \ref{asu:kert}, discussed in the next section.  
Let $\hat v_{1},\hat v_{2},...,\hat v_{n}$ denote the empirical functional principal components, EFPC's, of 
$Y_{1},Y_{2},...,Y_{n}$. Then, assuming the $Y_i$'s are centered, it is a basic fact of PCA that $\Span\{\hat v_1, 
\dots, \hat v_n\} = \Span\{Y_1,\dots, Y_n\}$. Recall that it is also a basic 
fact from linear algebra that the $\hat v_{1},\hat v_{2},...,\hat v_{n}$ can be 
completed to form a full orthonormal basis (all of the additional functions will 
have an empirical eigenvalue of 0). 
We then apply Parseval's identity to obtain
\begin{align*}
 &\sum_{i=1}^{n} \int_0^{ 1}\left (Y_{i}(t)- \int_0^{ 1} g(t,s,X_{i}(s)) ds 
\right)^2 dt 
 & = \sum_{i=1}^{n}\sum_{j=1}^{\infty}\left(\langle Y_{i},\hat v_{j}\rangle- 
\left \langle g, \int_0^{ 1}\int_0^{ 1} k_{t,s, X(s)} \hat v_{j}(t) dt ds 
\right \rangle_{\mbK}\right)^{2}.
\end{align*}
Define the subspace (of $\mbK$)
\[
\mathcal{H}_{1}= \Span \left\{ \int_0^{ 1}\int_0^{ 1} k_{t,s, 
X_{i}(s)} \hat v_{j}(t) dt ds , \hspace {0.1cm} i=1,2,...,n, \hspace 
{0.2cm} j= 1, \dots, n \right\},
\]
as well as its orthogonal compliment $ \mathcal{H}^{\perp}_{1}$. The space 
$\mbK$ can be decomposed into the direct sum: $\mbK = 
\mathcal{H}_{1}\oplus \mathcal{H}^{\perp}_{1}$, which means that we can write 
any function $g \in \mbK$ as $ g = g_{1}+g^{\perp}_{1} $, with $ g_{1}\in 
\mathcal{H}_{1} $ and $ g^{\perp}_{1} \in \mathcal{H}^{\perp}_{1}$. Using this 
decomposition we have that, for $1 \leq j \leq n$,
 \begin{align}
 \left \langle g,\int_0^{ 1} \int_0^{ 1} k_{t,s, X_i(s)} \hat v_{j}(t) dt ds 
\right \rangle_{\mbK}&= \left \langle g_{1}, \int_0^{ 1}\int_0^{ 1} k_{t,s, 
X_i(s)} \hat v_{j}(t) dt ds \right \rangle_{\mbK}.\label{Eq:tt}
\end{align}
Since $\Vert g\Vert^2_{\mbK}=\Vert g_1\Vert^2_{\mbK}+\Vert g^\perp_1\Vert^2_{\mbK}$, it follows from 
\eqref{e:rss} and \eqref{Eq:tt} that $ \hat{g}\in \mathcal{H}_{1}$ and 
so has the form 
\[
\hat g(t,s,x)= \sum_{i=1}^{n} \sum_{j=1}^{n}\alpha_{ij} \int_0^{ 1} \int_0^{ 
1}k\left((t,s,x);(t',s', X_{i}(s'))\right)\hat v_{j}(t') dt' ds'.
\]
Note that this same expression would hold if we replaced the $\{v_j(t)\}$ with $\{Y_j(t)\}$ (since they span the same space), however, it would not hold for an arbitrary basis.  We use the FPCs for computational reasons as we discuss at the end of the section.   
To compute the estimate, $\hat g$, we only need 
to compute the coefficients $\{\alpha_{ij}\}$. As usual, the coefficients 
$\alpha_{ij}$ can be computed via a type of ridge regression. Note that
\begin{align*}
& \left \langle \hat{g}, \int_0^{ 1}\int_0^{ 1} k_{t,s, X_i(s)} \hat v_{j}(t) 
dt ds \right \rangle_{\mbK} \\
 & = \sum_{i'=1}^{n} \sum_{j'=1}^{n}\alpha_{i'j'} \int_0^{ 1} \int_0^{ 1} 
\int_0^{ 1} \int_0^{ 1} \langle k_{t,s, X_{i'}(s)}, k_{t',s', X_{i}(s')} 
\rangle_{\mbK} \hat v_{j'}(t) \hat v_{j}(t') dt ds dt' ds' \\
& = \sum_{i'=1}^{n} \sum_{j'=1}^{n}\alpha_{i'j'} \int_0^{ 1} \int_0^{ 1} 
\int_0^{ 1} \int_0^{ 1} k(t,s, X_{i'}(s); t',s', X_{i}(s')) \hat v_{j'}(t) \hat 
v_{j}(t') dt ds dt' ds' .
\end{align*}
Define
\[
A_{iji'j'} = \int_0^{ 1} \int_0^{ 1} \int_0^{ 1} \int_0^{ 1} k(t,s, X_{i'}(s); 
t',s', X_{i}(s')) \hat v_{j'}(t) \hat v_{j}(t') dt ds dt' ds'.
\]
Turning to the norm in the penalty we can use the same arguments to show that
\[
\| \hat g\|_\mbK^2= \langle \hat g, \hat g \rangle_\mbK = \sum_{iji'j'} 
\alpha_{ij} A_{iji'j'} \alpha_{i'j'}.
\]
Thus the minimization problem can be phrased as
\[
\sum_{i=1}^n \sum_{j=1}^n \left(Y_{ij} - \sum_{i'j'} 
A_{iji'j'}\alpha_{i'j'}\right)^2 + \lambda 
\sum_{iji'j'}\alpha_{ij}A_{iji'j'}\alpha_{i'j'}.
\]
We now vectorize the problem by stacking the columns of $Y_{ij}$ and 
$\alpha_{ij}$, denoted as $\bY_V$ and $\balpha_V$. We also turn the array 
$A_{iji'j'}$ into a matrix $\bA_V$, by collapsing the corresponding dimensions. 
We can then phrase the minimization problem as
\[
(\bY_V - \bA_V \balpha_V)^\top (\bY_V - \bA_V \balpha_V) + \lambda 
\balpha_V^\top \bA_V \balpha_V.
\]
Thus, the final estimate can be expressed as
\[
\hat \balpha_V = (\bA_V^\top \bA_V + \lambda \bA_V)^{-1} \bA_V \bY_V.
\]

Note that we are estimating $n^2$ parameters and inverting an $n^2 \times n^2$ 
matrix. Thus for computational convenience, it is often useful to truncate 
the EFPCs at some value $J < n$. However, even without truncating this approach still has the potential to lead to less parameters than the basis methods of \citet{scheipl2015functional}, where the number of parameters to estimate is $m^3$, with $m$ being the number of basis functions used in their tensor product basis. In contrast, our approach yields $n^2$ parameters, and combined with an FPCA, this can be reduced to $nJ$ with relatively little loss in practical predictive performance.  There is also the possibility of using an eigen-expansion on $k$ 
to reduce the computational complexity even further \citep{flame}, though we don't pursue that here.  

\subsection{Alternative Domains} \label{s:pre}
While our work is focused primarily on the ``classic" function-on-function paradigm, we briefly mention in this section an easy way to modify the kernels to allow for more complex domains. In particular, one major concern brought up by a referee is when both $X_i(t)$ and $Y_i(t)$ are observed concurrently.  In that case, the classic approach would actually use future values of the covariate to predict present values of the outcome.  Interestingly, we need only make a very slight adjustment to the kernels to handle such a setting.

The goal here is to adjust the model such that
\begin{align}
Y_i(t) = \int_0^t g(t,s,X_i(s))\,ds + \vep_i(t) \qquad 0 \leq t \leq 1, \label{e:mod}
\end{align}
or equivalently to require that $g(t,s,X_i(s)) = 0$ if $s>t$.  More generally, we can allow the domain of $X$ used to predict $Y$ to change arbitrarily with $t$.  Let $\{A_t \subset [0,1]: 0\leq t \leq 1\}$ be a collection of (measurable) subsets of the unit interval.  Fitting \eqref{e:mod} is equivalent to taking $A_t = [0,t]$, which is what we use to highlight this approach in Section \ref{s:app}.  We aim to fit the more general model
\[
Y_i(t) = \int_{A_t} g(t,s,X_i(s))\,ds + \vep_i(t) \qquad 0 \leq t \leq 1.
\]  Interestingly, this can be done through a simple modification of the kernel.  In particular, we can define a new kernel as
\[
\tilde k(t,s,x, t', s', x') =  1_{s \in A_t} 1_{s' \in A_{t'}}k(t,s,x, t', s', x'). 
\]
A direct verification shows that $\tilde k$ is a valid reproducing kernel as long as the original $k$ was.  Then our estimate would take the form
\begin{align*}
\hat g(t,s,x) &= \sum_{i=1}^n \sum_{j=1}^n \int_0^1 \int_0^1 \tilde k(t,s,x;t',s',X_i(s')) \hat v_j(t') dt' ds' \\
& = 1_{s \in A_t} \sum_{i=1}^n \sum_{j=1}^n \int_0^1 \int_{s' \in A_{t'}}  k(t,s,x;t',s',X_i(s')) \hat v_j(t') ds' dt',
\end{align*}
which means that $\hat Y_{n+1}(t)$ can be computed using only $\{X_{n+1}(s):  s \in A_t\}$ and a very slight modification of our current approach.  We illustrate this technique in Section \ref{s:app}.

\section{Asymptotic Theory} \label{s:theory}

In this section, we demonstrate that the \textit{excess risk}, $\Re_n$ (defined below), of our 
estimator converges to zero at the optimal rate. Optimal convergence of 
$\Re_n$, for scalar-on-function linear regression was established by 
\citet{cai2012minimax}, while optimal convergence for the continuously 
additive scalar-on-function regression model was established in 
\citet{wang2015optimal}. In both cases an RKHS estimation framework was 
used. Because our model involves a functional response, the form of the excess 
risk $\Re_{n}$ is different and requires some serious mathematical extensions 
over previous works. However, we will show that the convergence rate for our 
model is the same as the one found in \citet{wang2015optimal}.

We begin by defining the excess risk, $\Re_n$. Let $X_{n+1}(t)$ be new 
predictor which is distributed as, but independent of $(X_i(t))^n_{i=1}$. We let $\E^*$ 
denote the expected value, conditioned on the data $\{(Y_i, X_i): 1 \leq i \leq 
n\}$. Then the excess risk is defined as
\[
\Re_n = \E^* \left[
\int^1_0 \int^1_0 (\hat g(t,s,X_{n+1}(s)) - g(t,s,X_{n+1}(s)))^2 \ dt ds
 \right].
\]
Note that $\Re_n$ is still a random variable as it is a function of the data. 
Intuitively, this quantity can be thought of as prediction error, namely, for a 
future observation, how far away is our prediction from the optimal one where 
the true $g$ is known.  For ease of exposition, we present all of assumptions below, even the ones discussed previously.   

\begin{asu} \label{a:main}
We make the following assumptions.  

\subasu \label{asu:model} The observations $\{Y_i(t), X_i(t)\}$ are assumed to satisfy
\[
Y_i(t) = \int g(t,s,X_i(s)) \ ds + \vep_i(t)
\]
where $\{X_i\}$ and $\{\vep_i\}$ are independent of each other and iid across $i=1,\dots,n$.

\subasu \label{asu:kernel} Denote by $\mcL_k$ the integral operator with $k$ as its kernel:
\[
(\mcL_k f)(t,s,x) := \int k(t,s,x; t', s', x') f(t',s',x') \ dt' ds' dx'.
\]
  The kernel, $k$, which also defines the RKHS, $\mbK$, is assumed to be symmetric, positive definite, and square integrable.

\subasu \label{asu:kert}Assume that there exists a constant $c>0$ such that for any $f \in \mbK$ and $t \in [0,1]$ we have 	
\[
 \E \left( \int_0^{ 1} f(t,s,X(s)) \  ds 
\right)^{4} \leq c \left[ \E\left( \int_0^{ 1} f(t,s,X(s)) \ ds\right)^{2}\right]^{2} < \infty.
\]

\subasu \label{asu:eigen} Let $\mcL_k^{1/2}$ denote a square--root of $\mcL$ (which exists due to Assumption \ref{asu:kernel}) 
and define $k^{1/2}_{t,s,x} := \mcL_k^{-1/2} k_{t,s,x}$.  Define the operator, $C$, as
\[
C(f) = \E \left[ \int \int \int
k^{1/2}_{t,s,X_i(s)} \langle k^{1/2}_{t,s',X_i(s')} , f\rangle_{L^2} \ ds ds' dt
\right].
\]
Assume that the eigenvalues $\{ \rho_{k}: k\geq 1 \}$ of  $C$ 
satisfy $ \rho_{k} \asymp k^{-2r}$ for some constant $ r>1/2$. 

\subasu \label{asu:var} There exists a constant $M>0$ such that, for all $t \in [0,1]$ and $i=1,\dots,M$
\[ 
 \E(\epsilon_i^2(t)) \leq M < \infty.
\]

\subasu \label{asu:ball}The function $g$ lies in $\Omega$, which we assume is a closed bounded 
ball in $\mbK$. 

\end{asu}

We are now in a position to state our main result.
\begin{theorem} \label{t:main}
If Assumption \ref{a:main} holds and the penalty parameter, $\lambda$, is chosen 
such that $\lambda \asymp n^{-\frac{2r}{2r+1}}$ then we have that
\[
 \lim_{A \to \infty} \lim_{n \to \infty} \sup_{g\in \Omega } \mathbb P \left 
(\Re_{n} \geq A n^{-\frac{-2r}{2r+1}} \right)=0.
\]
\end{theorem}

Before interpreting this result, let us discuss each of the assumptions individually.  Assumption \ref{asu:model} explicitly defines the model we are considering.  Assumption \ref{asu:kernel} ensures that the kernel has a spectral decomposition via Mercer's theorem, which will be used extensively. 
Assumption \ref{asu:kert} is fairly typical in these sorts of asymptotics, assuming that the fourth moment is bounded by a constant times the square of the second.  Assumption \ref{asu:eigen} introduces a central quantity that is used extensively in the proofs.  While not immediately obvious, this assumption basically states how ``smooth" or ``regular" the function $g$ is, as $g$ must lie in $\mbK$, whose kernel contributes to $C$.  In such results it is common for $X$ to contribute to the asymptotic behavior as the prediction error depends on the complexity of the $X$.  
Note that $k^{1/2}_{t,s,x}$ is a well defined quantity and it is easy to show via the reproducing property that it is an element of $L^2([0,1]^2 \times \mbR)$. 
The operator $C$ does depend on the choice of the square-root $\mcL_k^{1/2}$ (which is not a unique choice), however its eigenvalues do not.
Assumption \ref{asu:var} simply assumes that the point-wise variance of the errors is bounded, while the last assumption requires that the true function lie in a ball in $\mbK$, which is used to control the bias of the estimate.   
 

The rate given in Theorem \ref{t:main} is the same as was found in the scalar outcome case in \citet{wang2015optimal}, thus we know that this is the minimax rate of convergence.  In our case, as well as in \citet{wang2015optimal} and \citet{cai2012minimax}, it is the interaction between the covariance of $X$ and the kernel $k$ which determines the optimal rate.  The proof is quite extensive and given in the appendix.  The idea of the proof is to rephrase the estimate using operator notation instead of the representation theorem.  The difference between the estimate and truth is then split into a bias/variance decomposition.  Bounding the bias turns out to be relatively straight forward.  Bounding the variance is done by decomposing it into five more manageable pieces, and then bounding each of them separately.  Our task is complicated by the fact that the errors and response are now functions, where as in both  \citet{wang2015optimal} and \citet{cai2012minimax} they were scalars.  This requires extending many of the lemmas to this new setting, as well as using some completely new arguments to get the necessary bounds in place.


\section{Simulation Study}\label{s:sims}
Here we investigate the prediction performance of AFFR. We compare it with a linear model estimated in one of two ways.  The first way will be denoted as $LMR$ (linear model reduced) and $LMF$ (linear model full), where both use FPCA to reduce the dimension of the predictors, but $LMR$ also reduces the dimension of the outcome, while $LMF$ does not.  To implement our approach we relied heavily on the \textit{TensorA} package 
\cite{van2007tensora} in {\tt R}, which allowed us to carryout various tensor products very quickly.  

We consider three different settings for $g(t,s,x)$ one 
linear and two nonlinear forms:
\begin{enumerate} [(a)]
\item Scenario (a): $g(t,s,x) = tsx$,
\item Scenario (b): $g(t,s,x) = t+s+x^2$,
\item Scenario (c): $g(t,s,x) = tsx^2+x^4$.
\end{enumerate} 
 In all settings, the predictors $X_{i}(t)$ and 
errors $\epsilon_{i}(t)$ are taken to be iid Gaussian processes with mean 0 and the following covariance function from the Mat\'ern family:
\begin{equation*}
 C(t,s)= \left( 1 + \frac{\sqrt 5 |t - s|}{\rho} + \frac{5 |t - s|^2}{3 \rho^2} 
 \right)
 \exp\left( -\frac{\sqrt 5 |t - s|}{\rho} \right), 
\end{equation*}
where $\rho=1/4$.
For the RKHS we considered both the Gaussian kernel
\begin{equation*}
\begin{split}
 k\left((x,y,z),(x^{'},y^{'},z^{'}) \right) & =e^{ -\delta \left[(x-x^{'})^2 
+(y-y^{'})^2+(z-z^{'})^2 \right] },
\end{split}
\end{equation*}
and exponential kernel
\begin{equation*}
\begin{split}
k\left((x,y,z),(x^{'},y^{'},z^{'}) \right) & =e^{ -\delta \left[|x-x^{'}| 
	+|y-y^{'}|+|z-z^{'}| \right] },
\end{split}
\end{equation*}
where $\delta$ is the range parameter.  We will examine the sensitivity of our approach to this parameter in Tables \ref{t:two} and \ref{t:exp}.
 All of the curves ($X_{i}(t)$, $Y_{i}(t)$, and $\vep_i(t)$) were simulated on a $M=50$ equispaced grid between $0$ and $1$. The data is approximated using $K=100$ B-splines.  We denote by $J_{X}$ and $J_{Y}$ the number of principal components of 
$X$ and $Y$ respectively.  These steps are carried out using the {\it Data2fd} and {\it pca.fd} functions in the {\it R} package {\it fda}.   Our approach uses an FPCA on $Y$ only, but the LMR approach uses the FPCs for both $X$ and $Y$ .  The common recommendations for choosing $J_Y$ is either to use some cutoff for explained variability (commonly 85\%) or to look for an elbow in the scree plot ($J_X$ can also be chosen the same way or using a model based criteria such as BIC) \citep{KRBook}. Using an 85\% cutoff here results in 3 FPCs for our simulations, though we also include 6 and 9 to show that our approach is not very sensitive to this choice as long as a large proportion of variability is explained.    However, one should note the trade offs when choosing $J_Y$.  In general, the major gain in choosing a smaller $J_Y$ is faster computation, which is nontrivial for this problem.  The major loss is that one ``gives up" on some proportion of the variability in $Y$.  For example, if the FPCs explain 95\% of the variability, then one immediately gives up on predicting that remaining 5\%.  This is a different consideration than when choosing FPCs for predictors.  In general, users can tailor this choice to their data; if one expects very accurate predictions then a larger $J_Y$ can be helpful so that one does not lose prediction accuracy, while if it is known a-priori that the prediction accuracy will be low, then $J_Y$ can be safely made smaller.

To evaluate the different approaches, we used 1000 repetitions of every scenario.  In each case we generate 150 curves to fit the different models and then generated another 150 curves to evaluate \textit{out-of-sample} prediction error.  The metric for determining prediction performance we denote as RPE, for relative prediction error.  This metric denotes the improvement of the predictions over just using the mean, and can be thought of as a type of out-of-sample $R^2$.  An RPE of 0 implies that the model shows no improvement over just using the mean, while an RPE of 1 means the predictions are perfect.  More precisely, we first compute the Mean Squared Prediction error as:
\[
MSPE = \sum_{i=1}^n \| Y_i - \widehat Y_i\|_{L^2}^2,
\]
where $\widehat Y_i$ is a predicted value using one of the three discussed models or simply the mean.  The RPE is then defined as
\[
RPE = \frac{MSPE_{mean} - MSPE}{MSPE_{mean}},
\]
where $MSPE_{mean}$ denotes the MSPE using a mean only model.  Note that even in the mean only model, all parameters are estimated on the initial 150 curves and prediction is then evaluated on the second 150.  Therefore, it is actually possible to have a numerically negative RPE if an approach isn't predicting any better than just using the mean.
 
The RPEs of $LMR$ and $LMF$ for 
the three models (a), (b), and (c) are summarized in Table \ref{t:one}. For both models, we took $J_X = 3$, which explained over $85\%$ of the variability of the predictors and for $LMR$ we took $J_Y=3$ PCs for the outcome as well.  
The RPEs for our approach with $\delta=\{2^{-3},2^{-2},2^{-1},1,2 
\}$ and $J_Y = 3, 6, 9$ are summarized in Tables \ref{t:two} and \ref{t:exp}, which represent the Gaussian and exponential kernels respectively.
An initial look at the tables confirms much of what one would expect.  When the true model is linear, the two linear approaches work best, resulting in about twice the RPE of AFFR.  However, when moving to the two nonlinear models, the AFFR approach does substantially better.  This increased performance is seen for any choice of $J_Y$ and $\delta$.  Furthermore, the prediction performance seems relatively robust to the choice of $J_Y$, $\delta$, and even the kernel.  In the case of $J_Y$ this is not so surprising as over $90\%$ of the variability of the $Y_i$ is explained by the first three FPCs.  In contrast, there is some sensitivity to the choice of $\delta$, but it is relatively weak given how much we are changing $\delta$ in each row.  In our application section we set $\delta$ using a type of median, but one could also refit the model with a few different $\delta$ and choose the one with the best prediction performance.  Given how consistent the AFFR predictions are, trying a few $\delta$ appears to be satisfactory, and large grid searches can be avoided.    

\begin{table}[ht]
\centering
\small
\begin{tabular}{|l|l|l|l|r|r|r|}
 \hline
 & Scenario (a) & Scenario (b)& Scenario (c)\\ 
 \hline
$ LMR$ & 0.045 & 0.030 & 0.060\\ 
$LMF$ & 0.045 &0.029 & 0.060\\ 
 \hline
\end{tabular}
\caption{Relative prediction errors, RPE, for the two linear models.  For both, the number of FPCs for the predictor is $J_X=3$.  LMR also reduces the dimension of the outcome with $J_Y=3$ FPCs.
\label{t:one}}
\end{table}

\begin{table}[ht]
\scriptsize
\centering
\begin{tabular}{|l|l|l|l|l|l|l|l|l|l|}
\hline
 & \multicolumn{3}{c|}{Scenario (a)} & \multicolumn{3}{c|}{Scenario (b)} & 
\multicolumn{3}{c|}{Scenario (c)} \\
\cline{1-10}
 & $J_{Y}=3$ & $J_{Y}=6$ & $J_{Y}=9$ &$J_{Y}=3$ & $J_{Y}=6$ & $J_{Y}=9$ & 
$J_{Y}=3$ & $J_{Y}=6$ & $J_{Y}=9$\\

\hline
$\delta=2^{-3}$ & 0.025 & 0.026& 0.026& 0.379 & 0.379 & 0.379 &0.840 & 0.840 & 
0.845 \\
$\delta=2^{-2}$ & 0.024 & 0.025 & 0.025&0.370 & 0.370& 0.370&0.816 & 0.804 & 
0.815\\
$\delta=2^{-1}$ & 0.023 & 0.024 & 0.023&0.360 & 0.361& 0.361&0.847 & 0.831 & 
0.830\\
$\delta=2^{0}$ & 0.022 & 0.023 & 0.021&0.346& 0.347& 0.347&0.83& 0.83 & 0.83\\
$\delta=2^{1}$ & 0.020 & 0.021 & 0.019&0.328& 0.328& 0.400&0.808& 0.808 & 
0.790\\
$PEV$ & 90.45\% &99.12\% & 99.88\%&90.82\% &99.10\% &99.84\%&91.25\% & 99.22\% 
&99.87\%\\
\hline
\end{tabular}
\caption{ Relative prediction error, RPE, for AFFR using a Gaussian kernel and with different kernel parameter values, $\delta$.  In every case the penalty parameter, $\lambda$, is chosen using cross-validation.  PEV indicates the proportion of explained variance of $Y$ for the 
corresponding number of FPCs, $J_{Y}$.\label{t:two}}
\end{table}

\begin{table}[ht]
	\scriptsize
	\centering
	
	\begin{tabular}{|l|l|l|l|l|l|l|l|l|l|}
		
		\hline
		& \multicolumn{3}{c|}{Model (a)} & \multicolumn{3}{c|}{Model (b)} & \multicolumn{3}{c|}{Model (c)} \\
		\cline{1-10}
		& $J_{Y}=3$ & $J_{Y}=6$ & $J_{Y}=9$ &$J_{Y}=3$  & $J_{Y}=6$ & $J_{Y}=9$ & $J_{Y}=3$  & $J_{Y}=6$ & $J_{Y}=9$\\
		
		\hline
		
		
		$\delta=2^{-3}$ & 0.021  & 0.023& 0.020& 0.368 & 0.379 & 0.379 &0.774 & 0.789& 0.775 \\
		
		$\delta=2^{-2}$ &  0.022 & 0.023 & 0.023&0.361 & 0.357& 0.359&0.813 & 0.805 & 0.815\\
		
		$\delta=2^{-1}$ &  0.022 & 0.023 & 0.023&0.350 & 0.349& 0.351&0.829 & 0.813 & 0.818\\
		
		$\delta=2^{0}$ &  0.021 & 0.022 & 0.022&0.338& 0.332& 0.334&0.780& 0.800 & 0.792\\
		
		$\delta=2^{1}$ &  0.020 & 0.019 & 0.019&0.300& 0.304& 0.302&0.743& 0.752 & 0.749\\
		
		$PEV$ &  90.45\% &99.12\% & 99.88\%&90.82\% &99.10\% &99.84\%&91.25\% & 99.22\% &99.87\%\\
		
		\hline
	\end{tabular}
	\caption[Table of the $RPP^{(NAFFR)}$ for Models (a),(b) and (c).]{ Relative prediction error, RPE, for AFFR using an exponential kernel and with different kernel parameter values, $\delta$.  In every case the penalty parameter, $\lambda$, is chosen using cross-validation.  PEV indicates the proportion of explained variance of $Y$ for the 
		corresponding number of FPCs, $J_{Y}$.\label{t:exp}}
\end{table}

\begin{figure}[htb]
	\centering
	\includegraphics[width=0.242\textwidth]{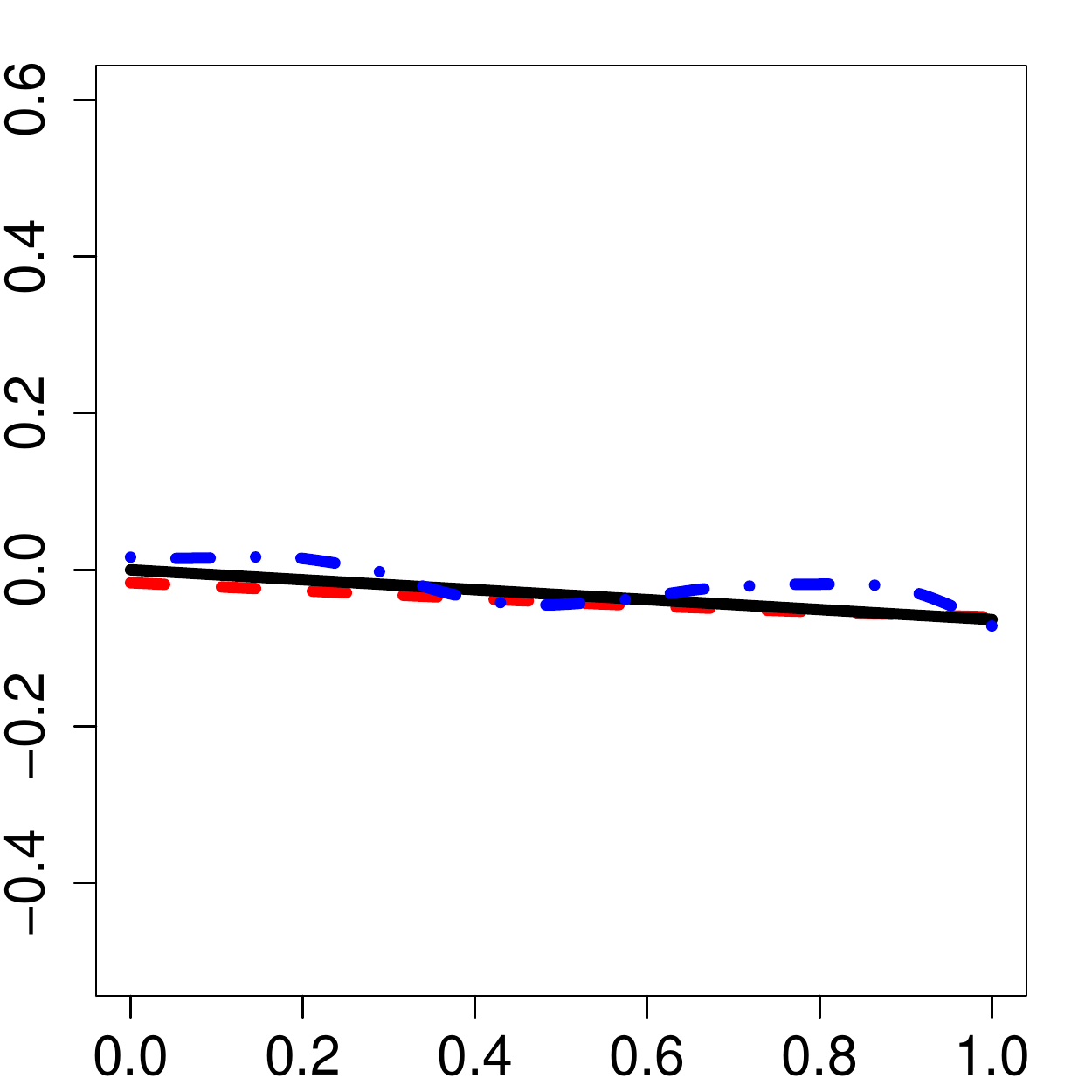}
	\includegraphics[width=0.242\textwidth]{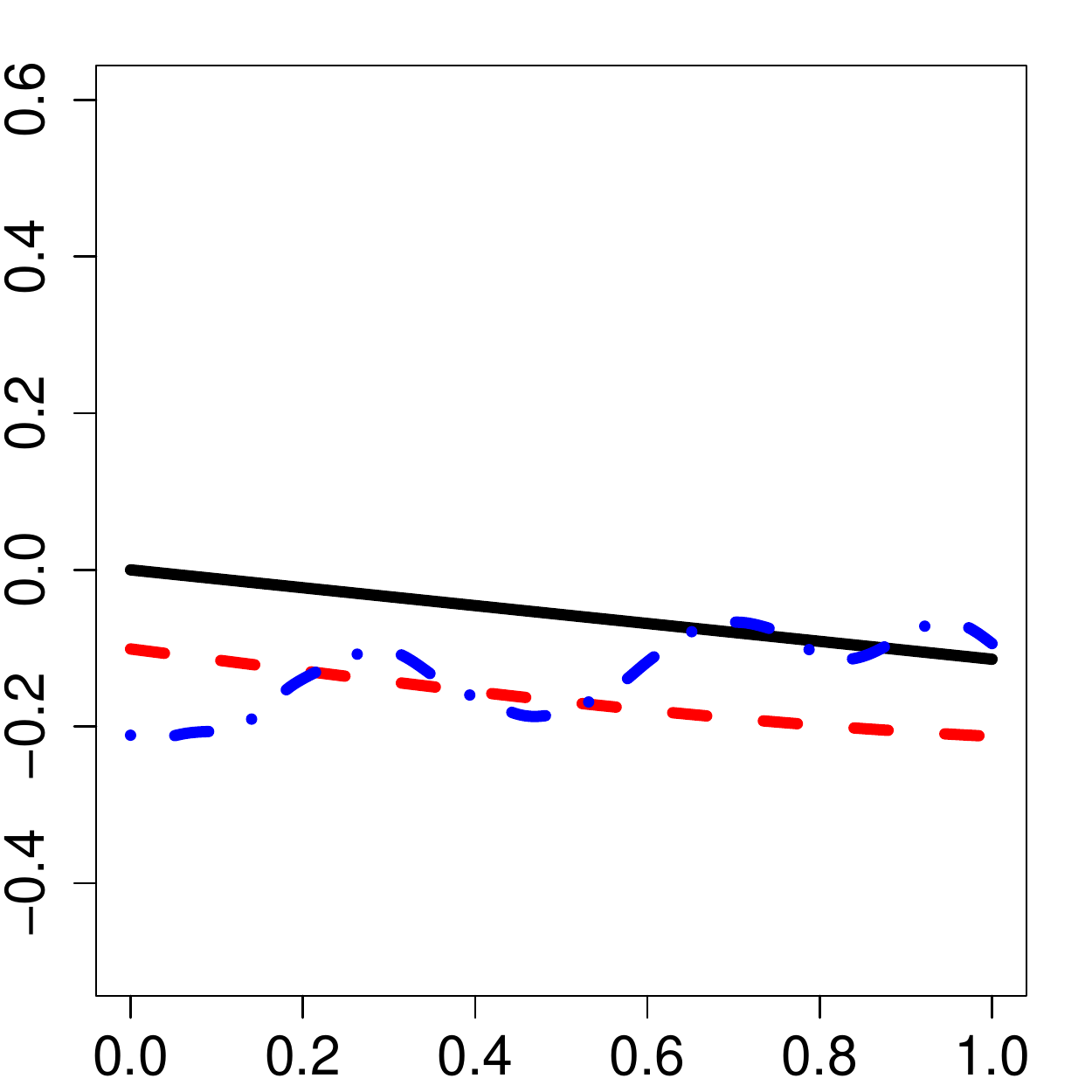}
	\includegraphics[width=0.242\textwidth]{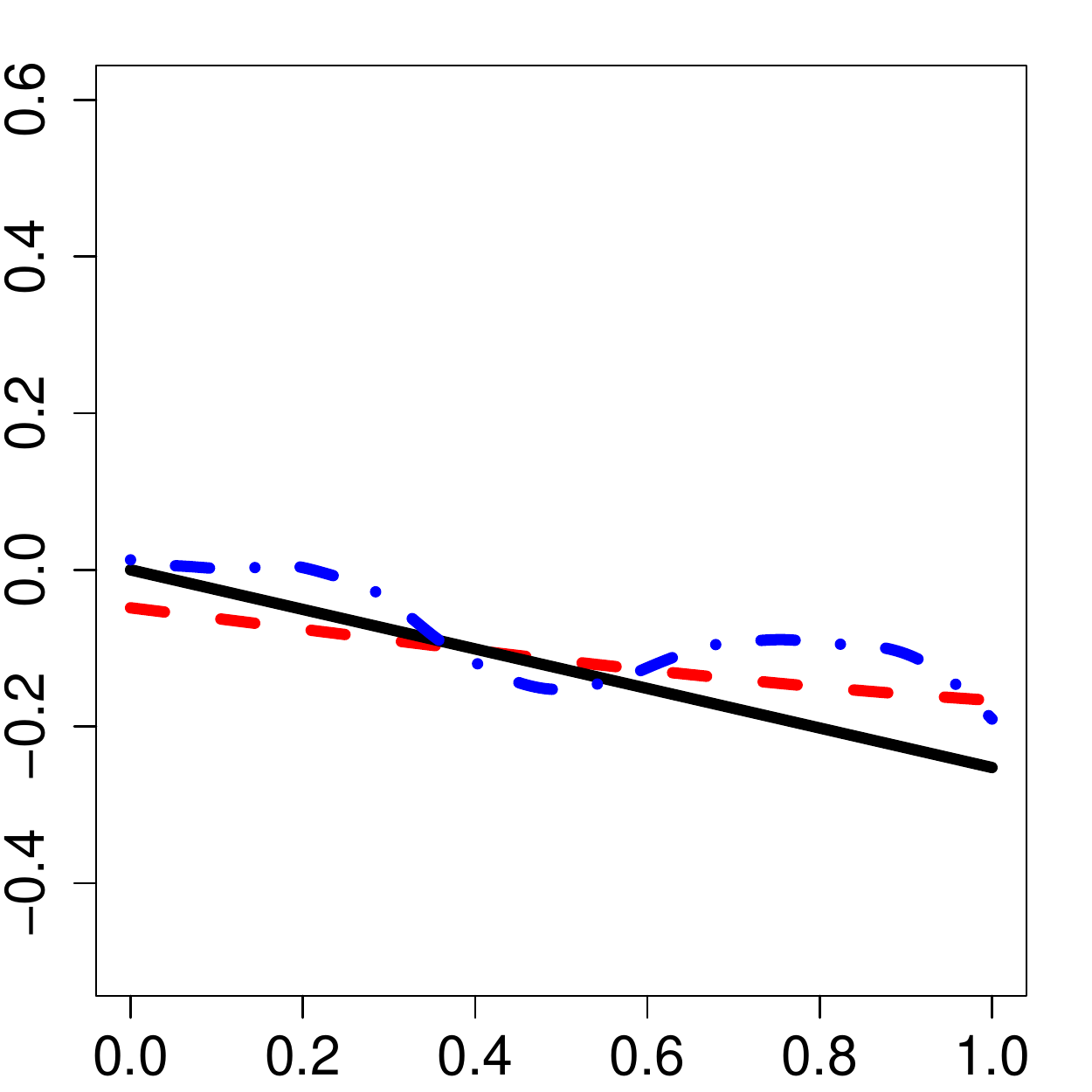}
	\includegraphics[width=0.242\textwidth]{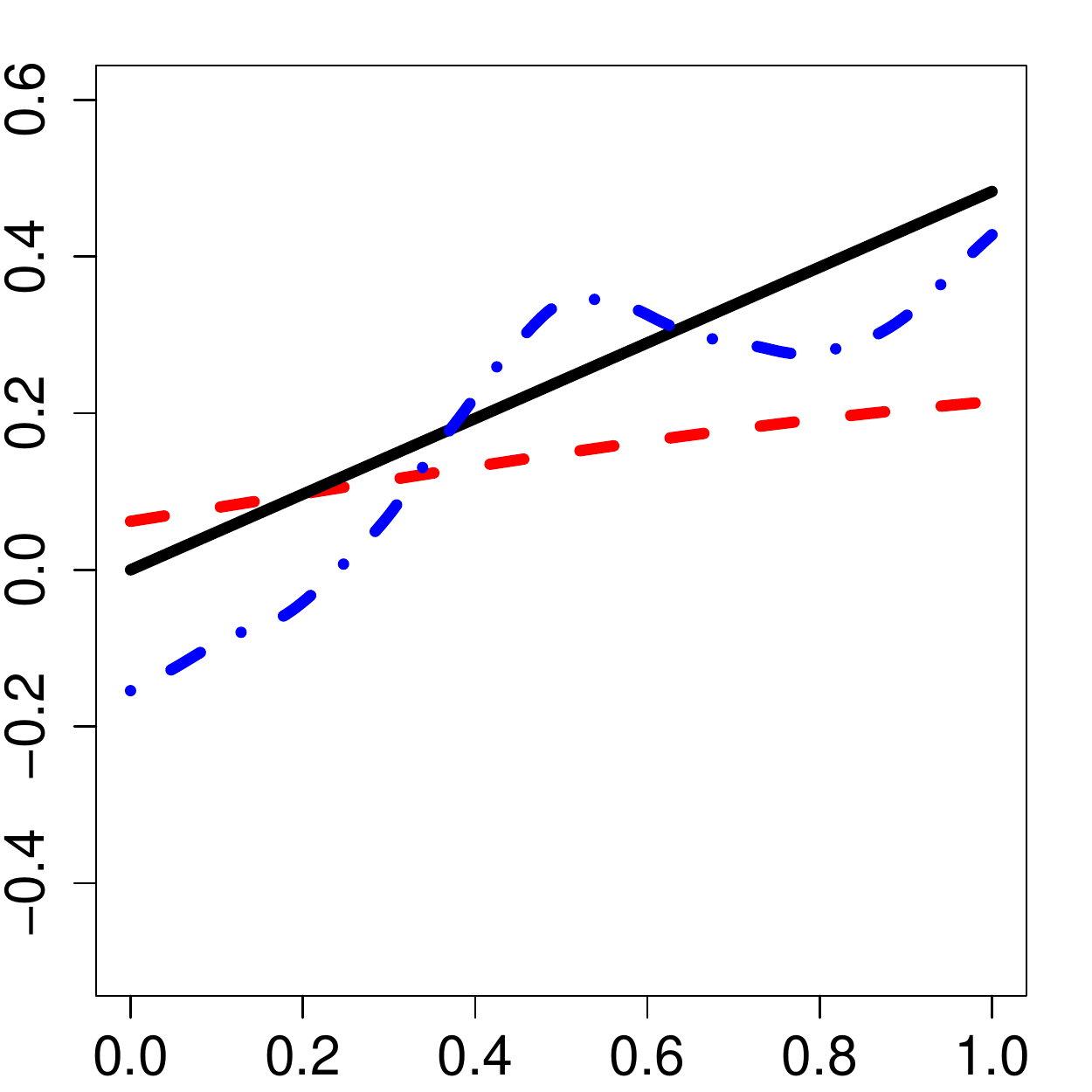}
	\includegraphics[width=0.242\textwidth]{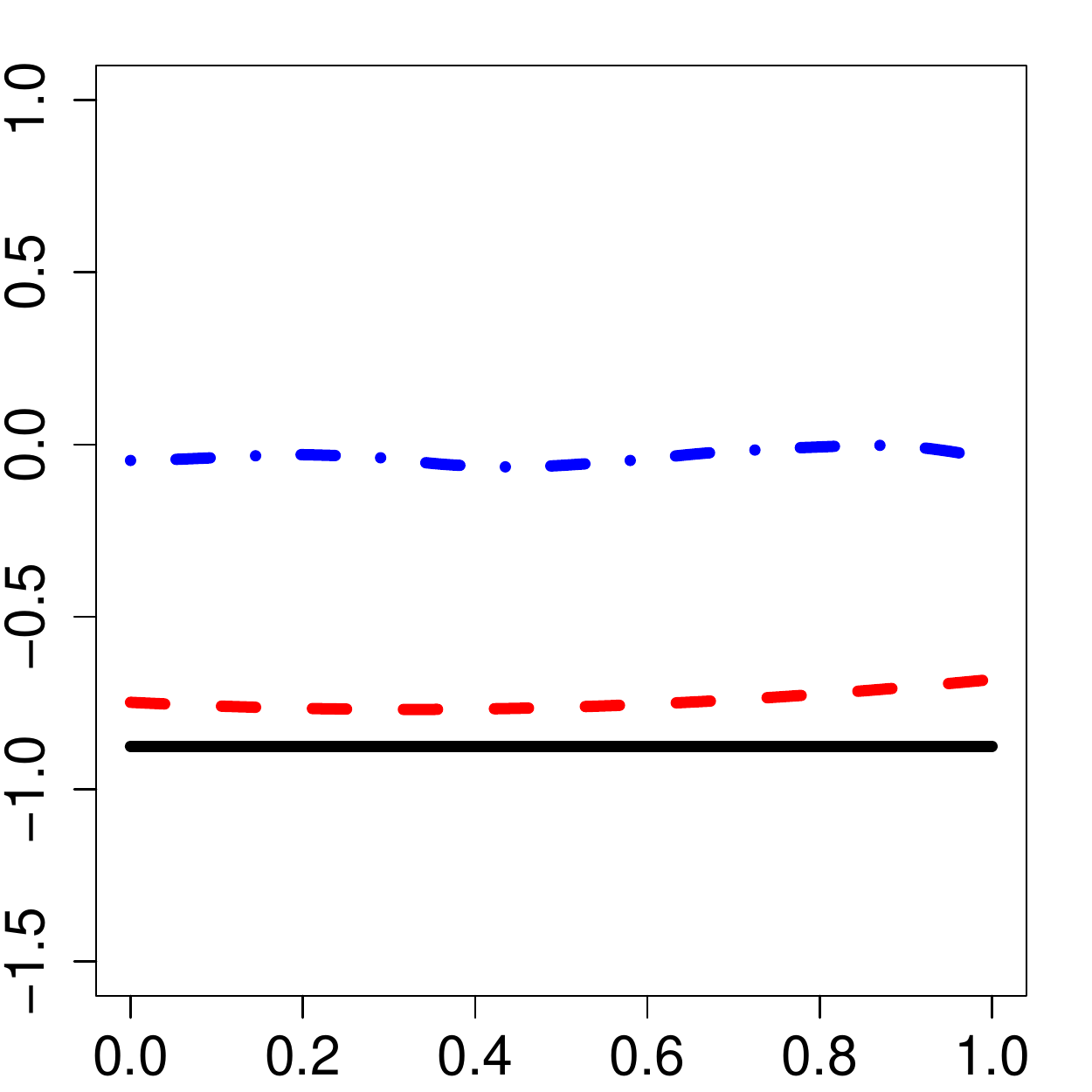}
	\includegraphics[width=0.242\textwidth]{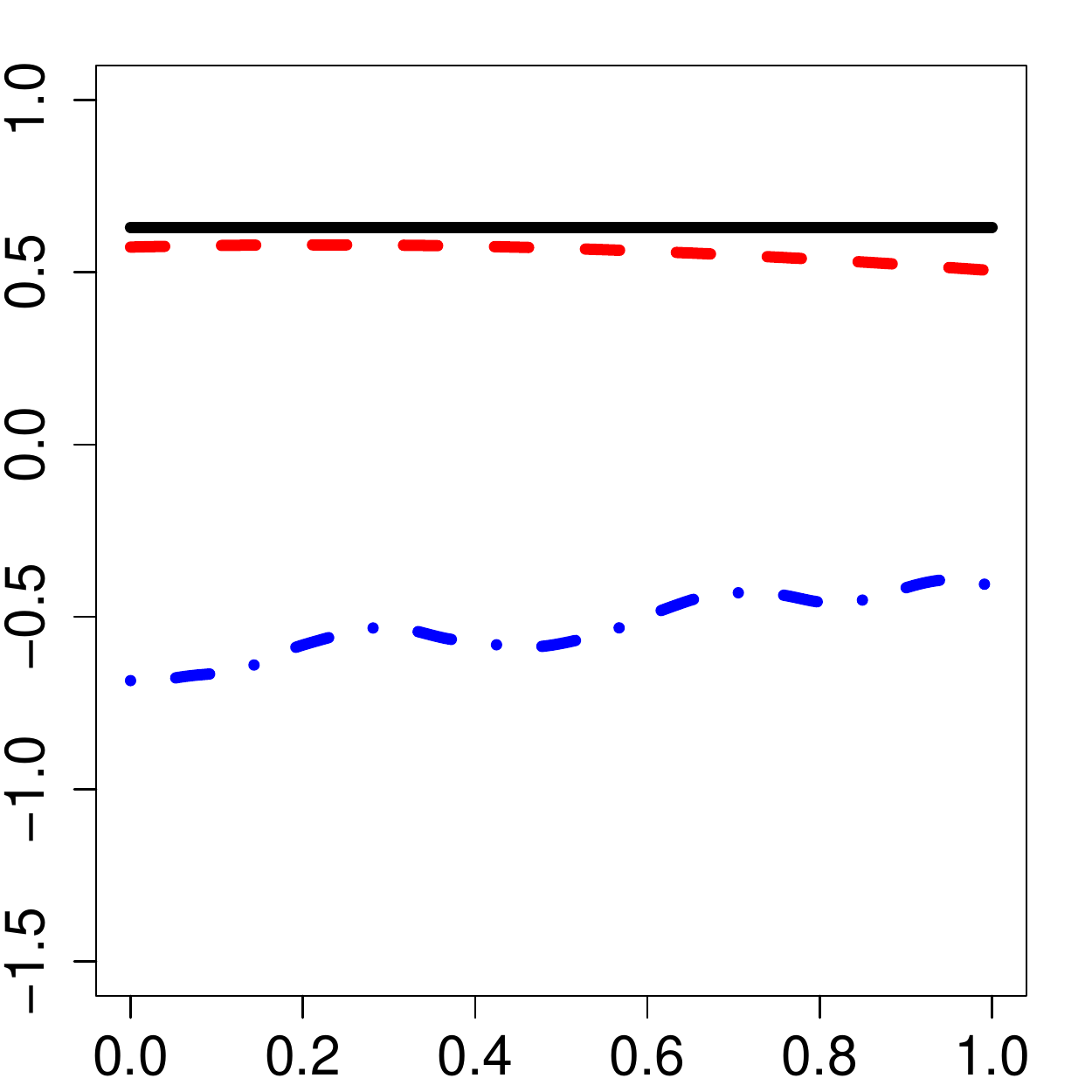}
	\includegraphics[width=0.242\textwidth]{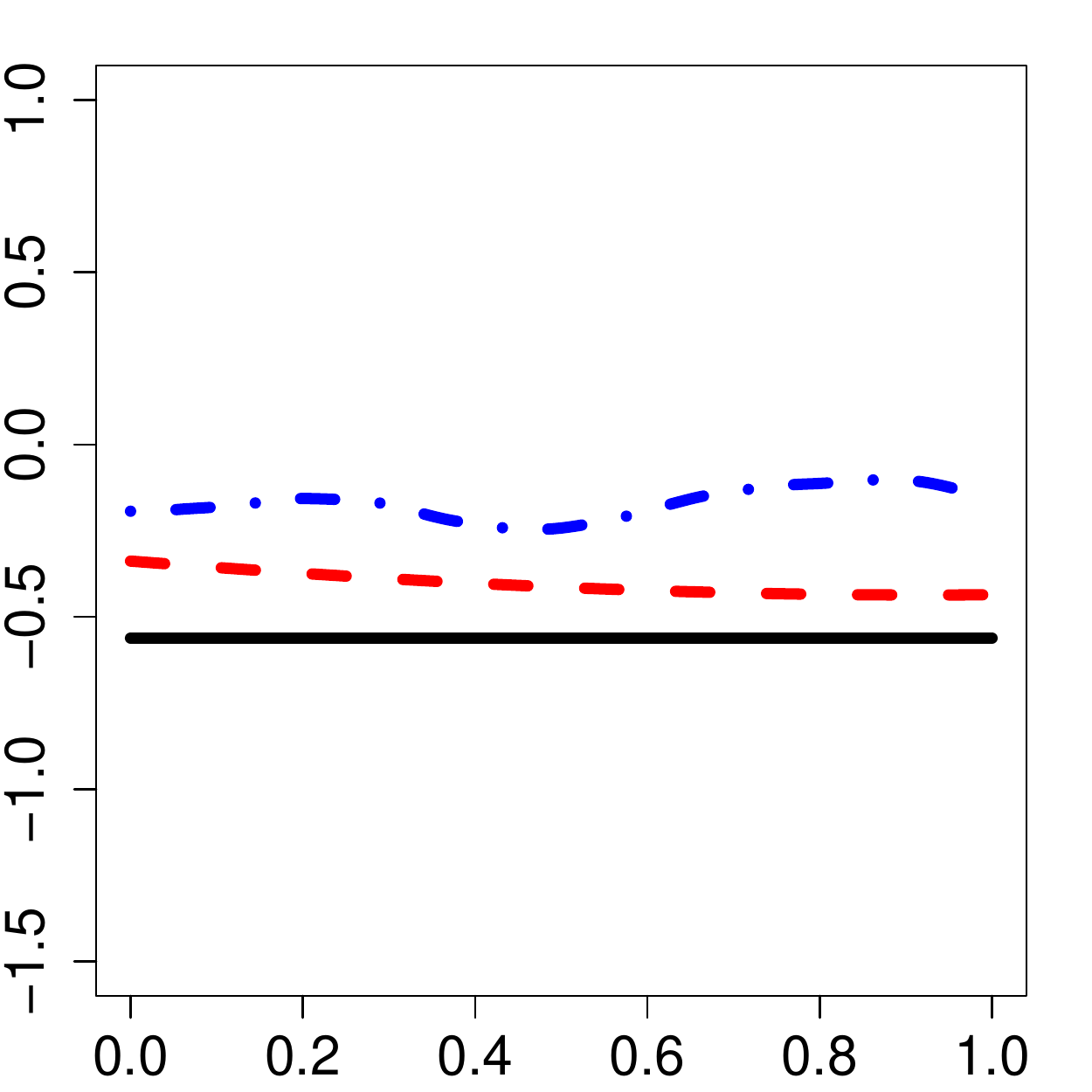}
	\includegraphics[width=0.242\textwidth]{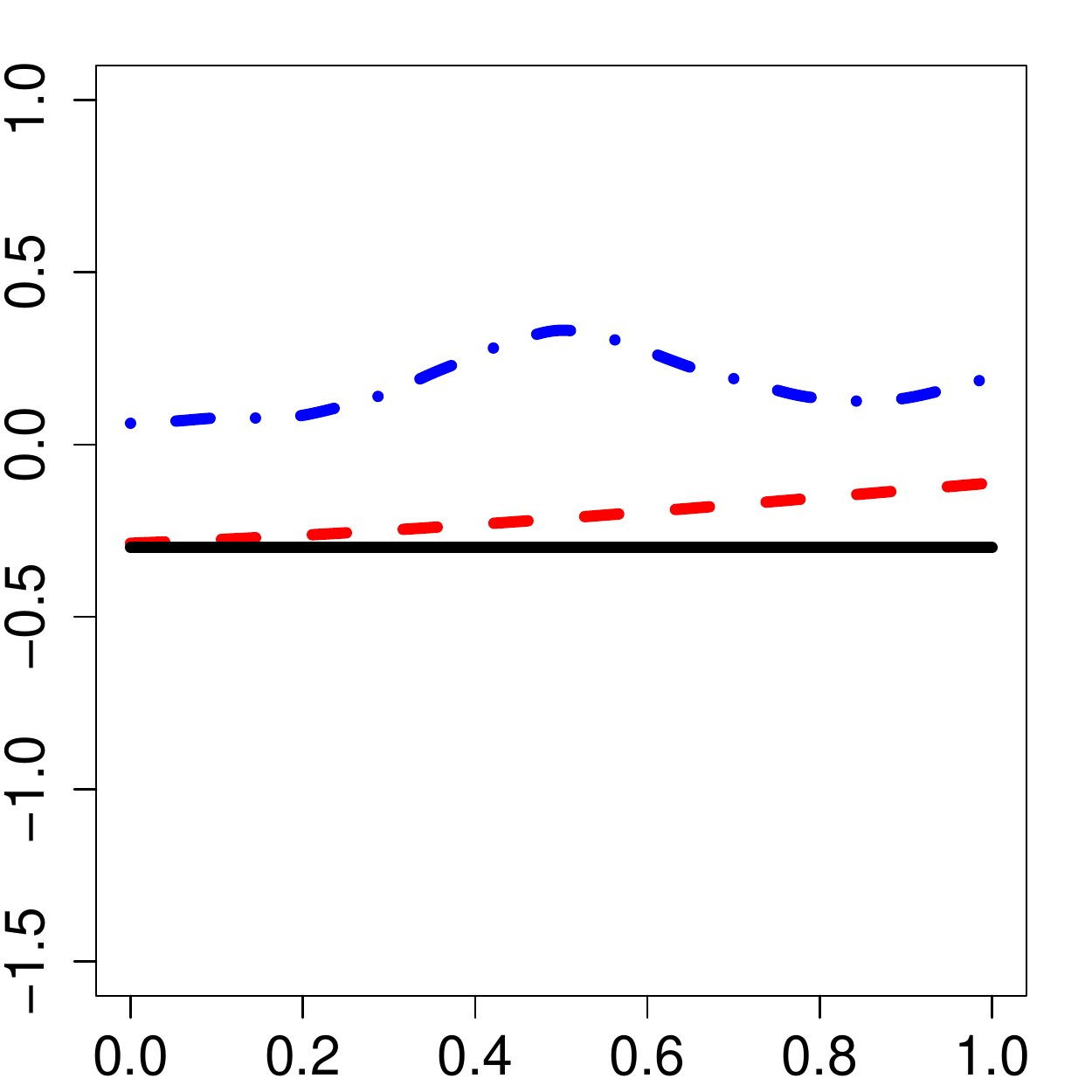}
	\includegraphics[width=0.242\textwidth]{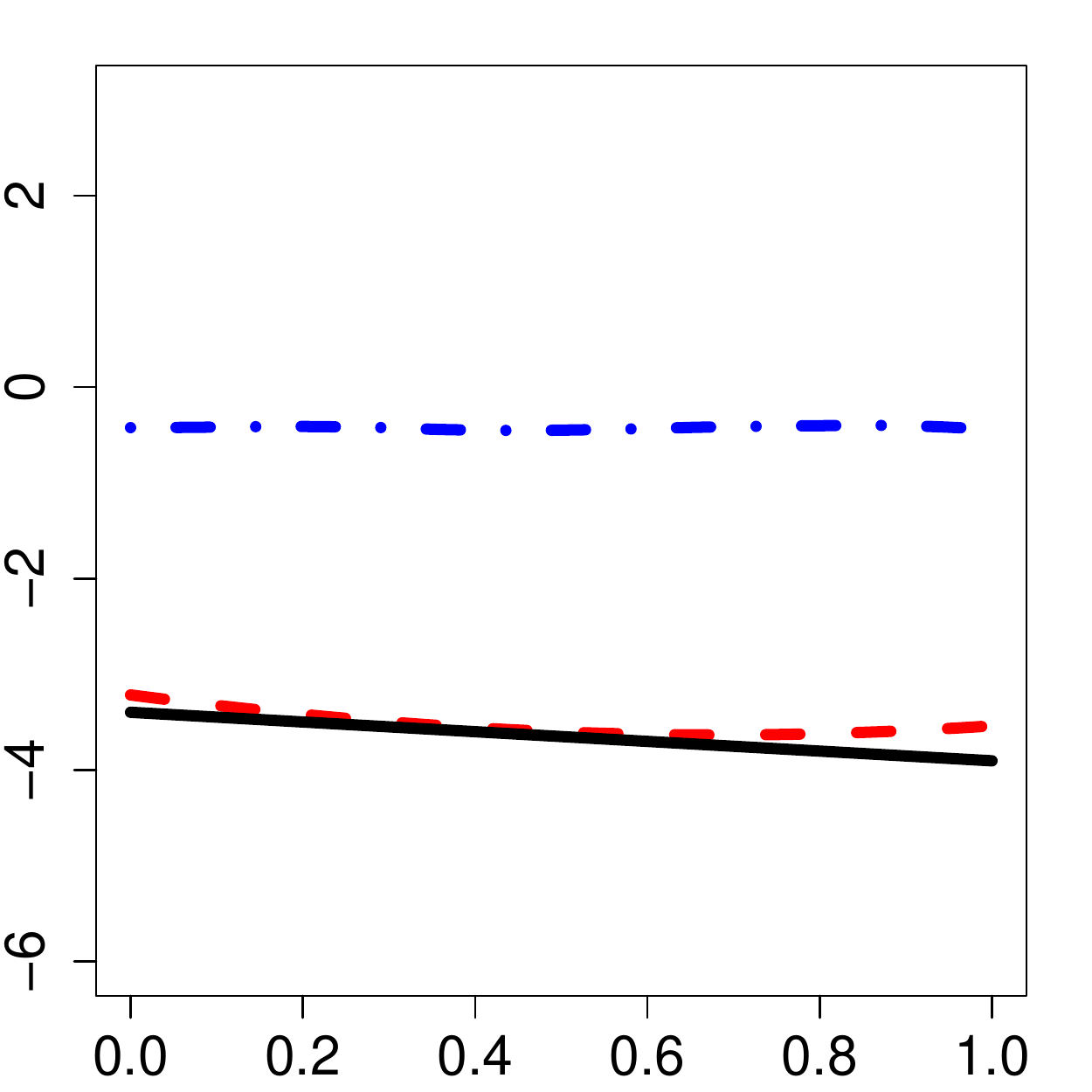}
	\includegraphics[width=0.242\textwidth]{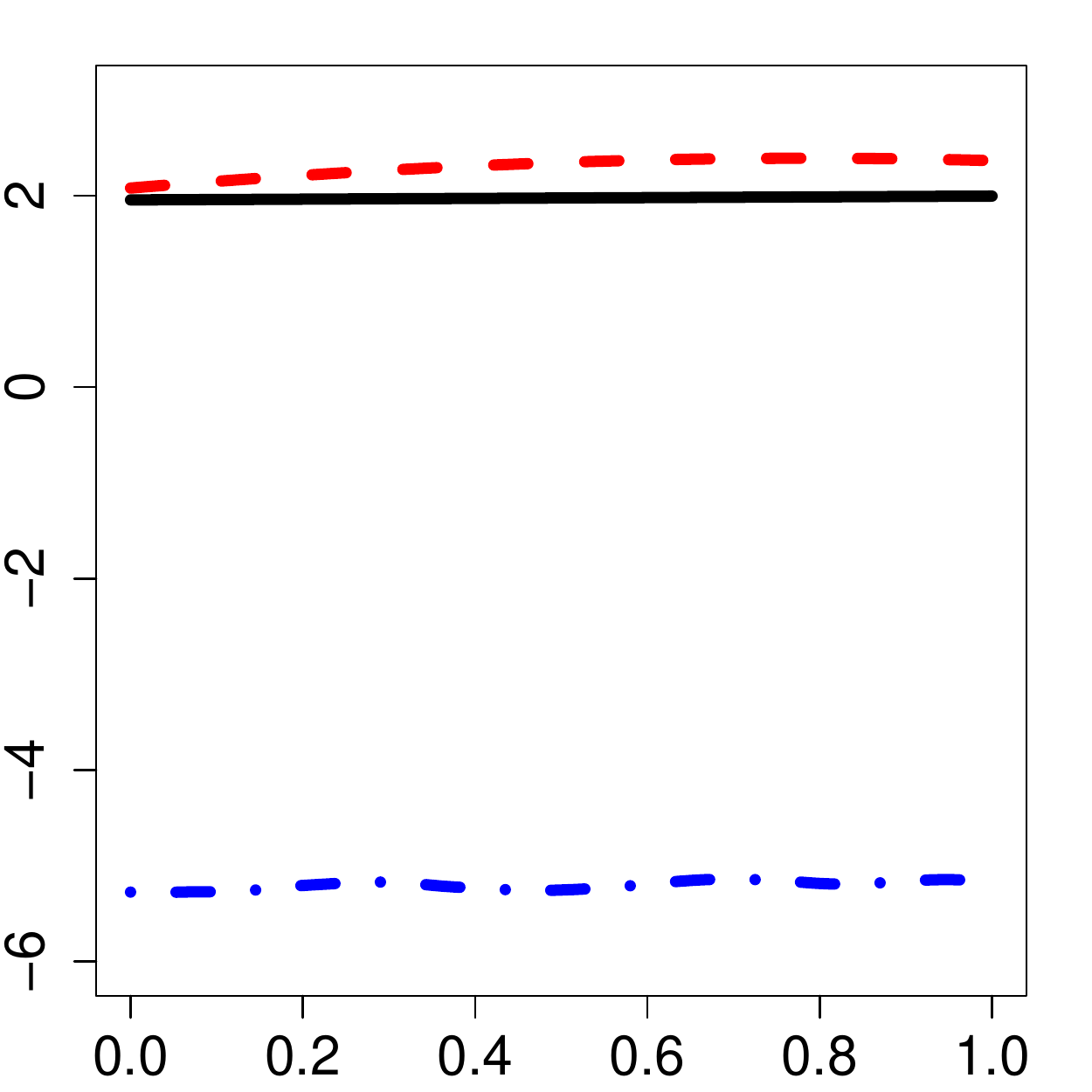}
	\includegraphics[width=0.242\textwidth]{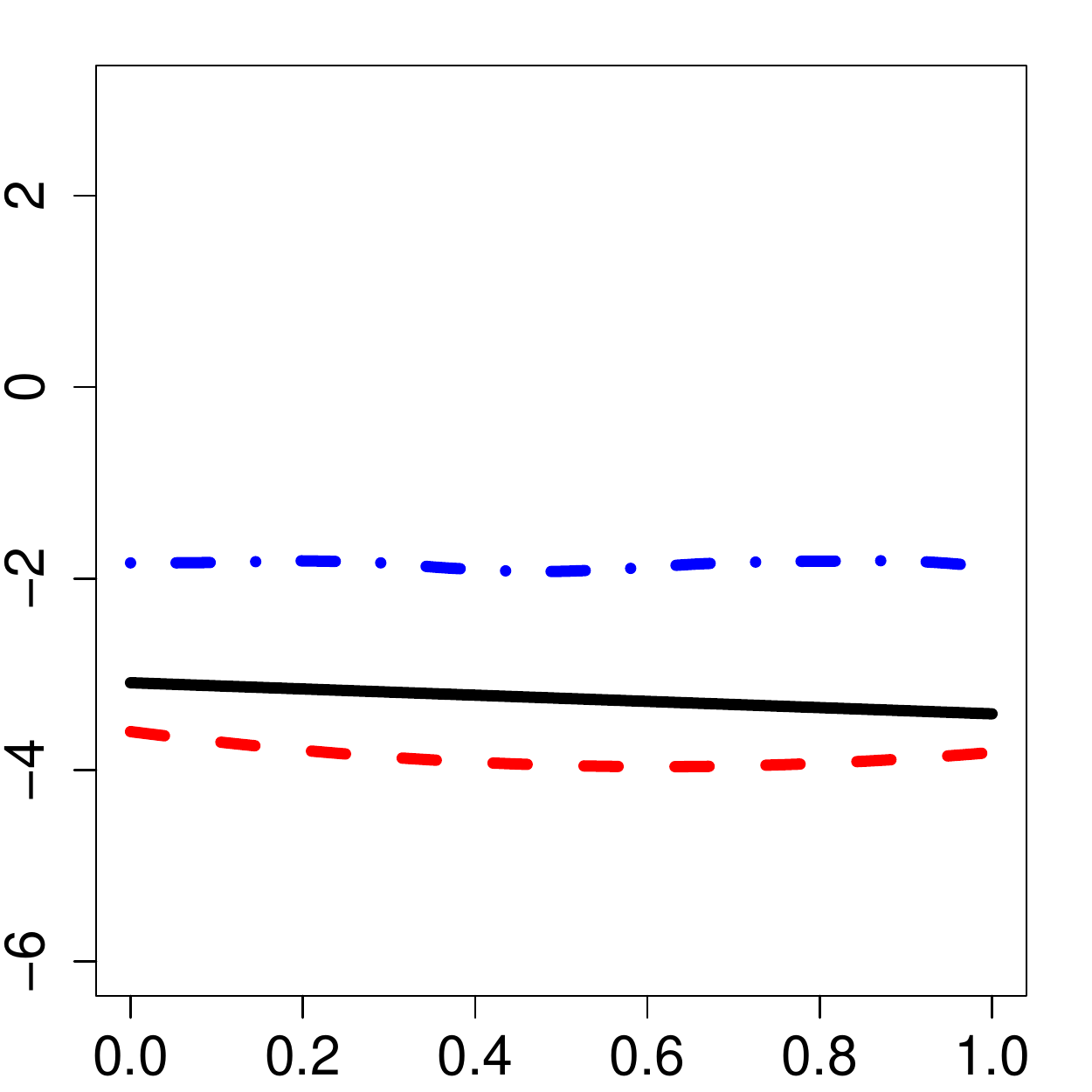}
	\includegraphics[width=0.242\textwidth]{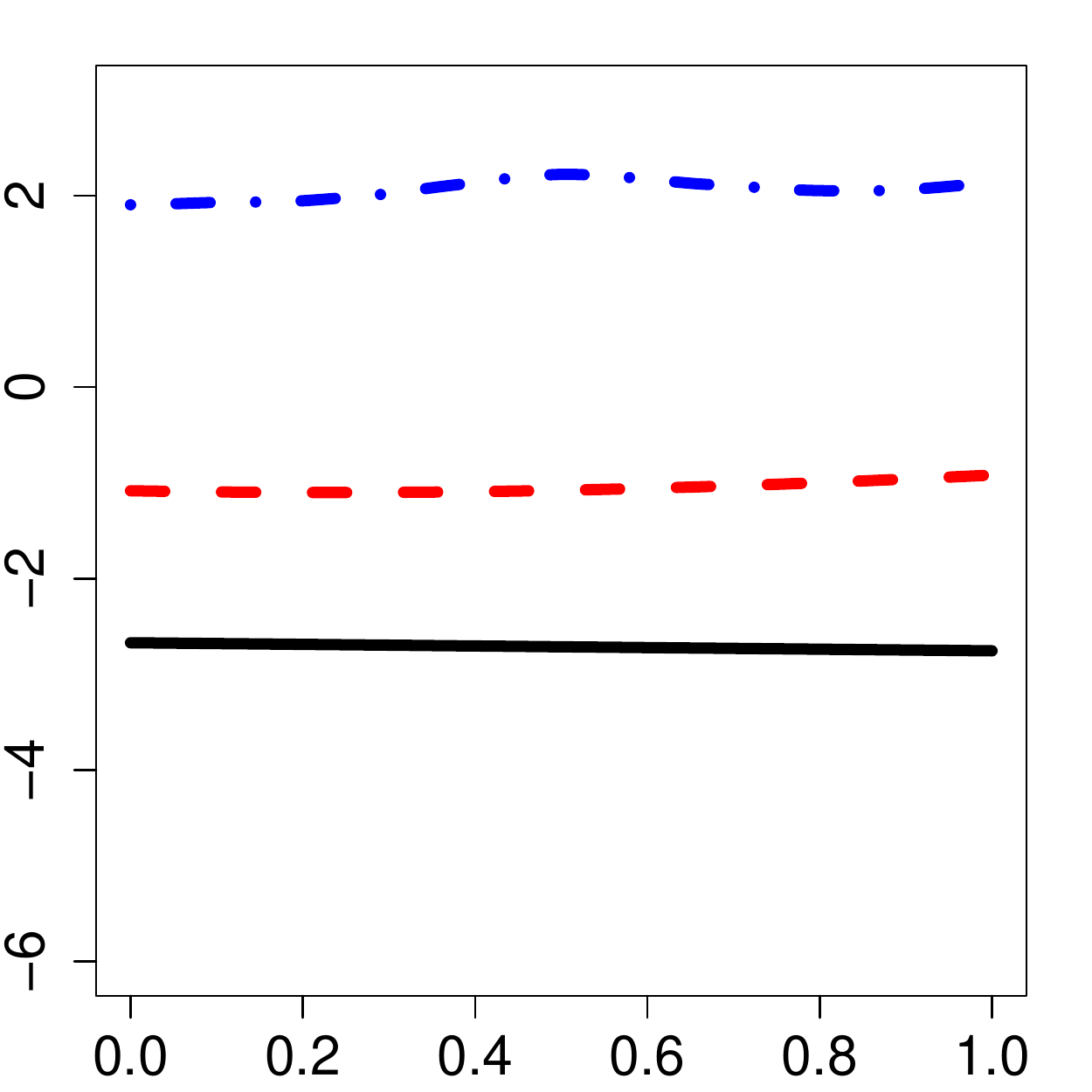}
	\caption{Plots of the optimal prediction  $E[Y(t)|X]$ (black), prediction using AFFR $\hat Y(t)$ (red dashed), and prediction using the unreduced linear model $\hat Y_{LM}(t)$ (blue dashed). The four plots on the top row correspond to the scenario (a), which is linear. The four plots in the middle row correspond to the scenario (b), which is nonlinear. The four plots in the bottom row correspond to the scenario (c) which is also nonlinear. \label{f:yhat}}
\end{figure}

\begin{figure}[htb]
	\includegraphics[width=0.327\textwidth]{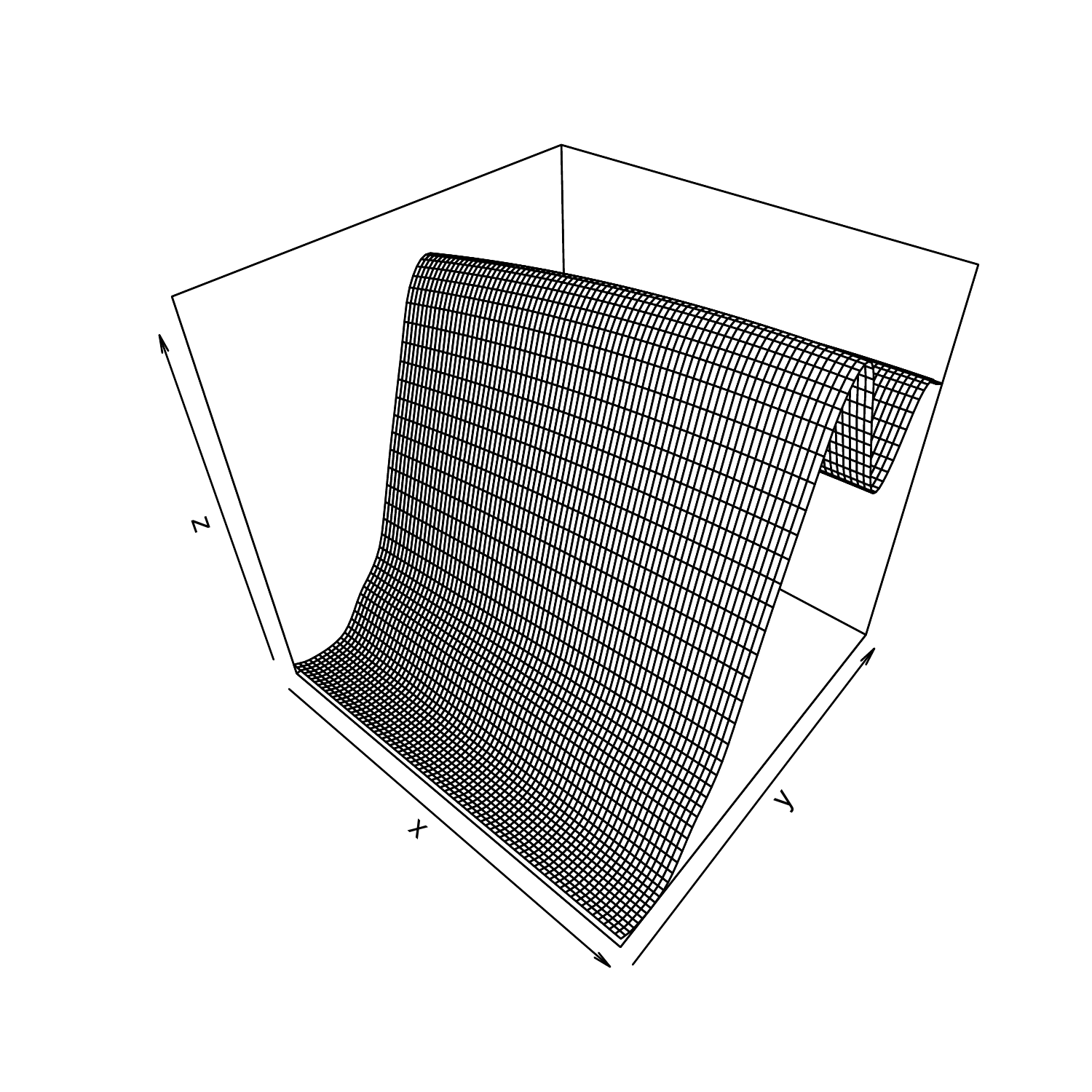}
	\includegraphics[width=0.327\textwidth]{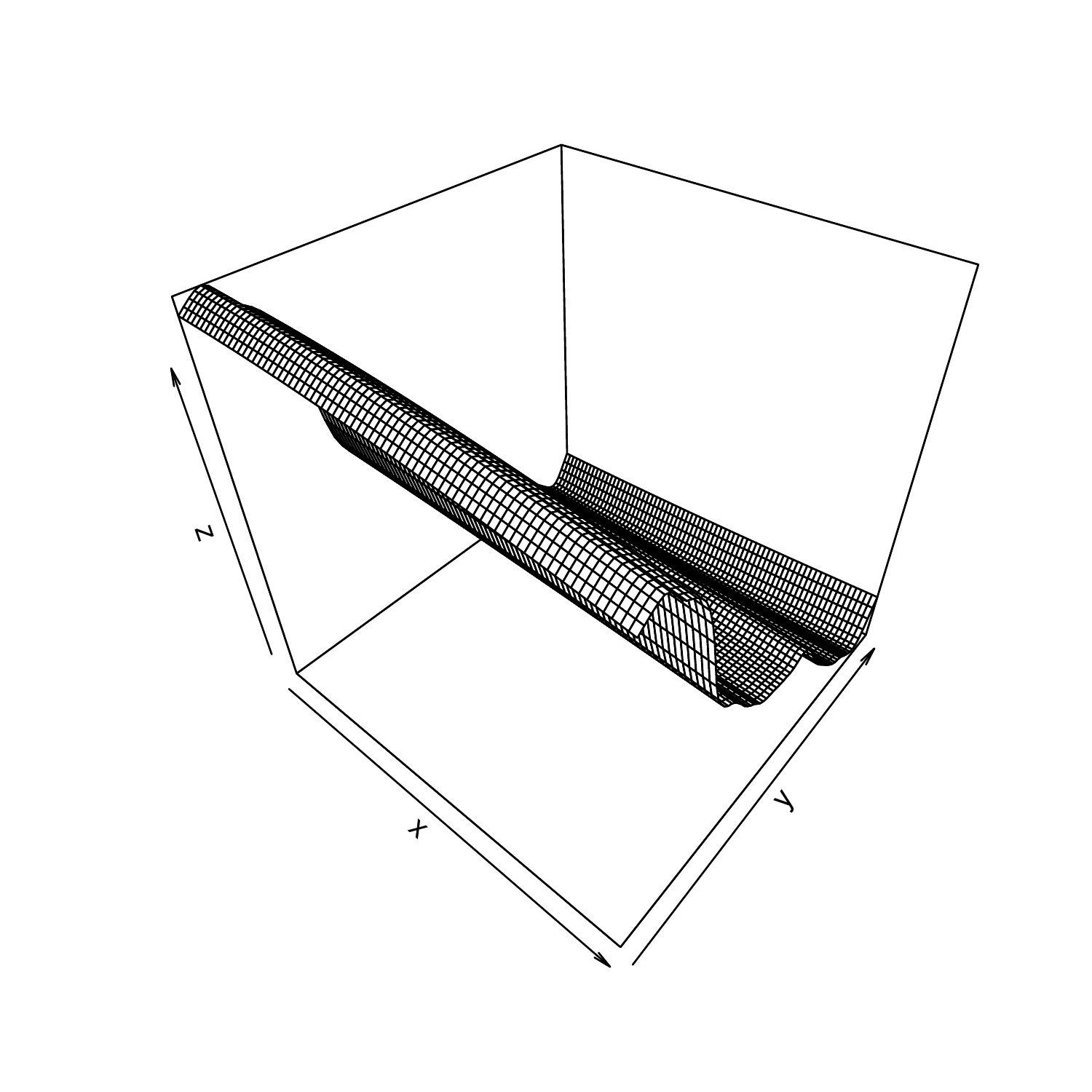}
	\includegraphics[width=0.327\textwidth]{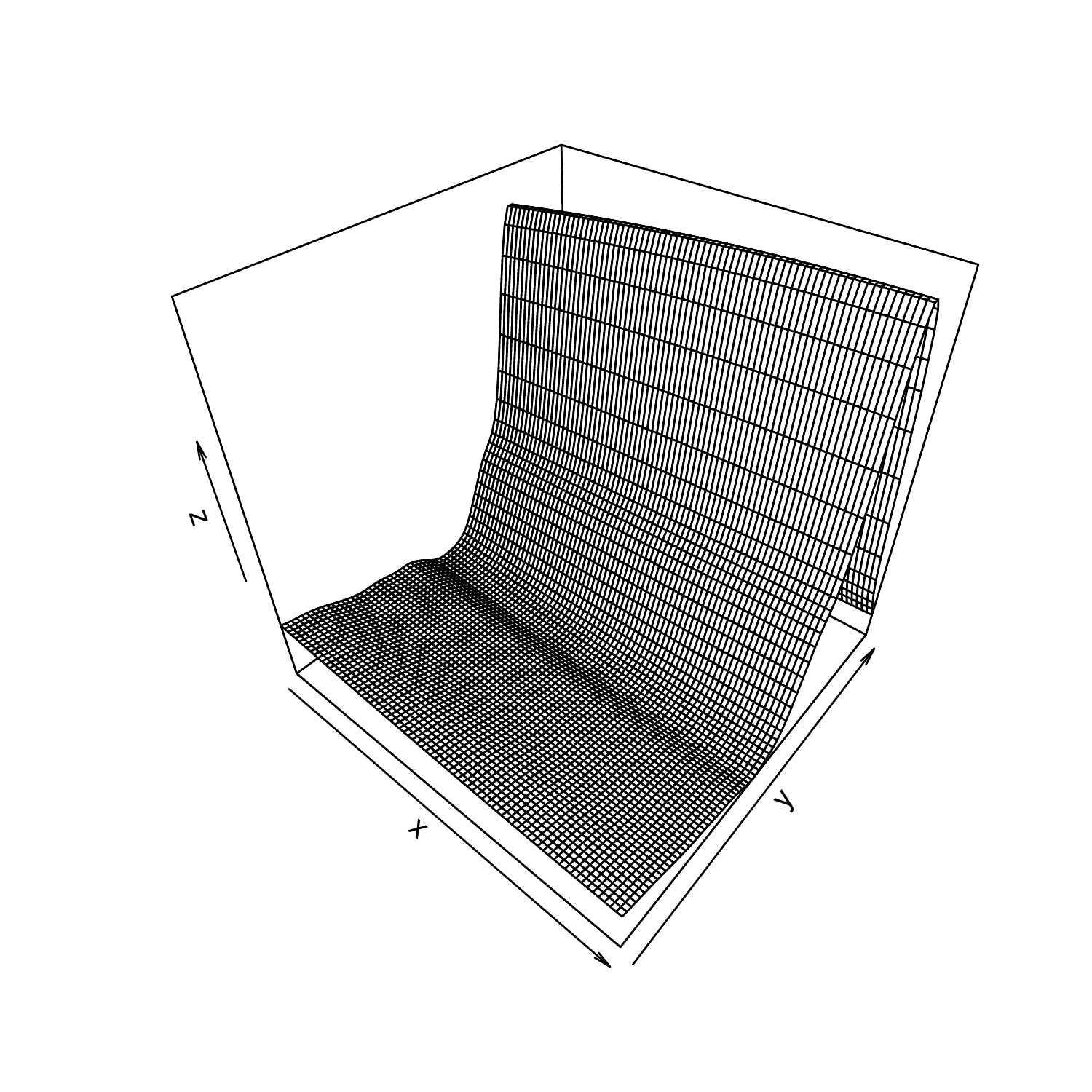}
	\includegraphics[width=0.327\textwidth]{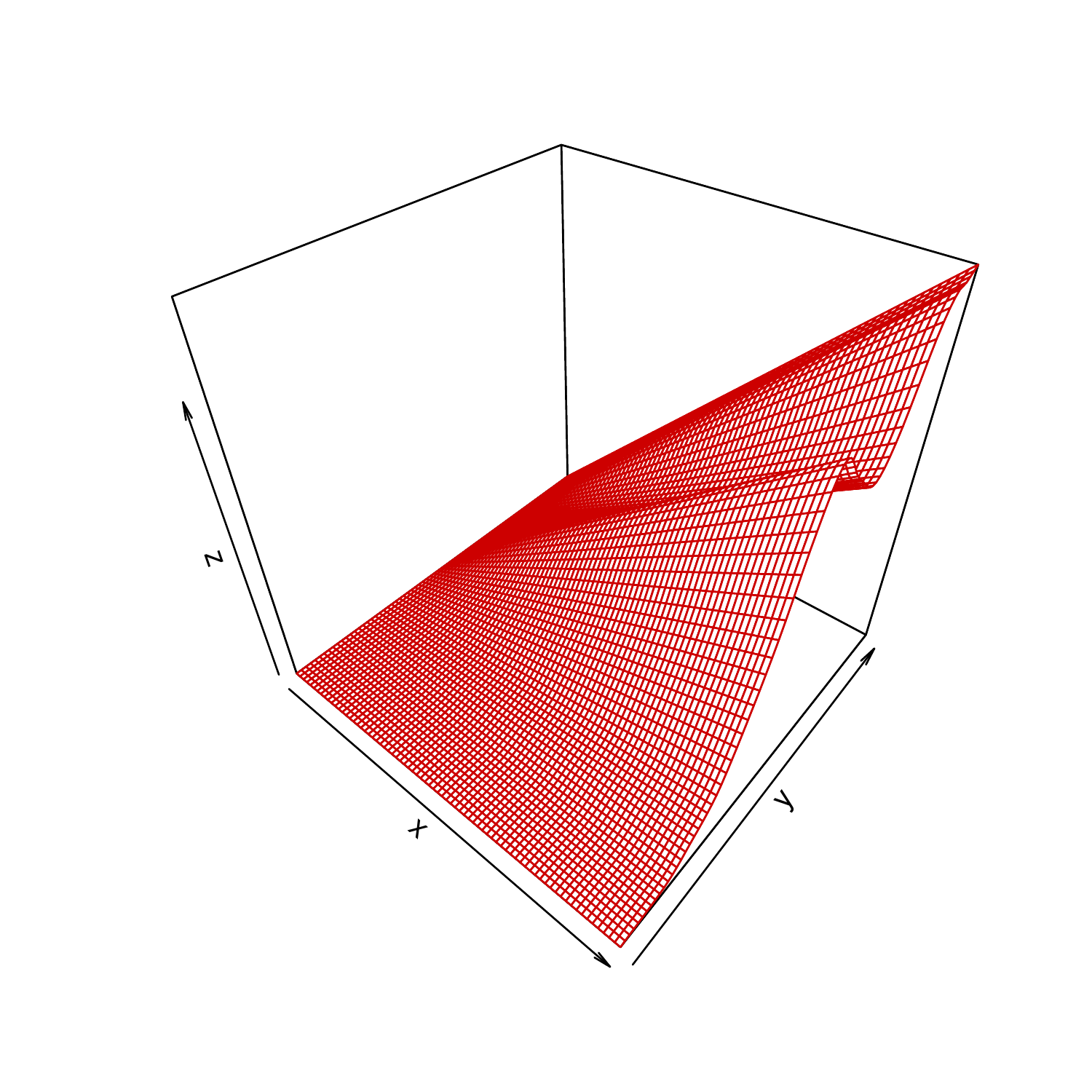}
	\includegraphics[width=0.327\textwidth]{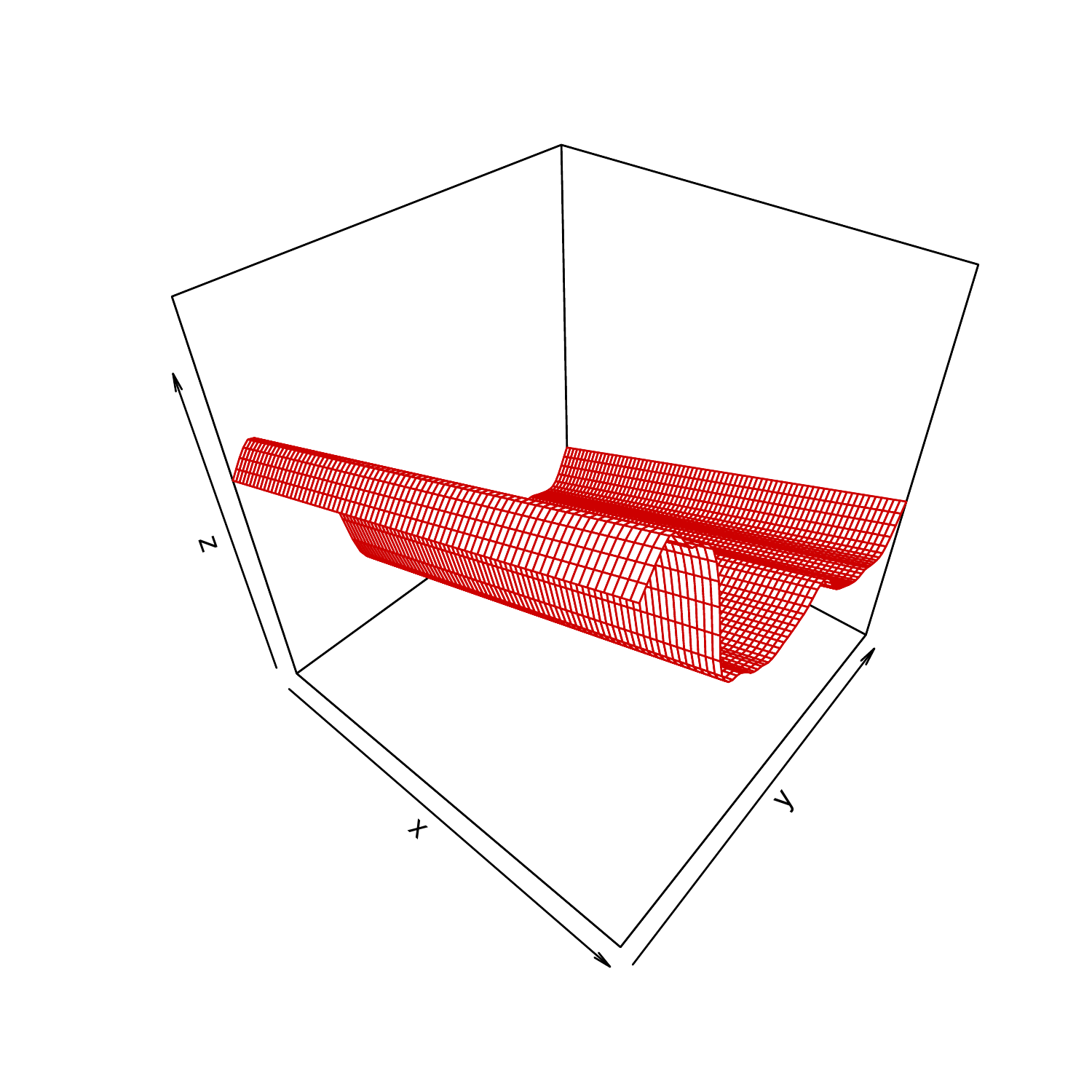}
	\includegraphics[width=0.327\textwidth]{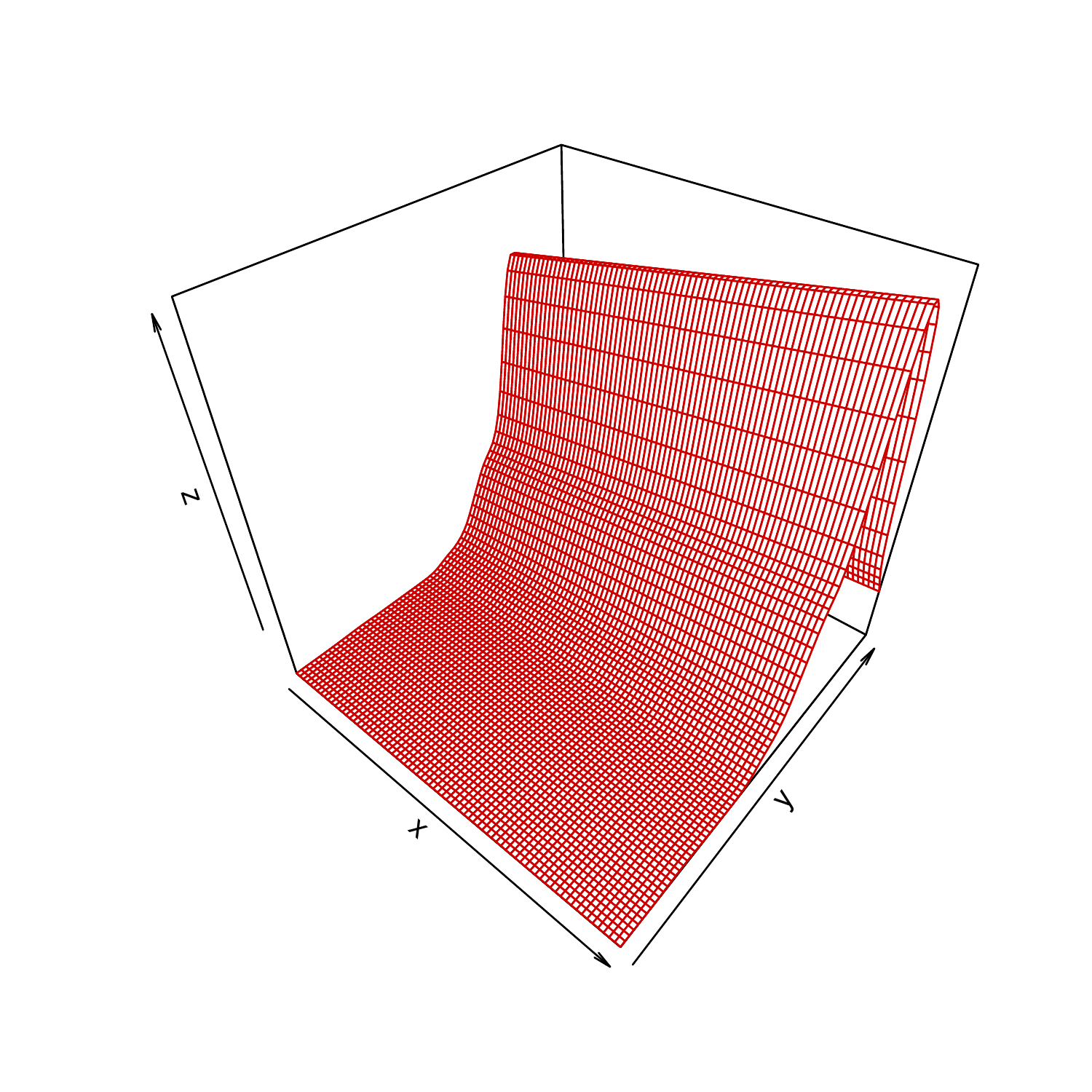}
	\caption{The top row plots one realization of $g(t,s,X(s))$ for models (a), (b), and (c) respectively.  The bottom row plots the corresponding (out of sample) prediction $\hat g(t,s,X(s))$. \label{f:ghat} }
\end{figure}

As a final illustration of the efficacy of AFFR, we provide several plots to help visualize the performance.  In Figure \ref{f:yhat} we plot several realizations of $Y_i$ and their corresponding (out of sample) predictions using the optimal prediction, $\E[Y(t) | X]$, AFFR, and the linear model without reducing the dimension of the $Y$.  We consider only the Gaussian kernel and take $\delta = 1/4$.  For the nonlinear scenarios (rows 2 and 3), one can clearly see the RPE results reflected in the predictions as AFFR is much closer to the optimal prediction.  In Figure \ref{f:ghat} we plot several realizations of $\hat g(t,s,X_i(s))$, which are again done out of sample along with the true value of $g(t,s,X_i(s))$.  Plotting in this way allows us to visualize $g$ using surfaces, where as plotting $g(t,s,x)$ would be challenging since the domain has three coordinates.  As we can see, the estimates are quite close to the true values, capturing the nonlinear structure quite well.

\section{Application to Cumulative Intraday Data}\label{s:app}
We conclude with an illustration of our approach applied to real data.  
 Cumulative Intra-Day Returns (CIDR's) consist of daily stock prices that are normalized to start at zero at the beginning of 
 each trading day.  FDA methods have been useful in analyzing such data \citep{gabrys2010tests,kokoszka2013predictability,horvath2014testing}, given the density at which stock prices can be observed.
Let $P_i(t_j)$ denote the price of a stock on day $i$ and time of day $t_j$.  The CIDRs are then defined as
\[
 R_{i}(t_{j}) = 100 \left[ \ln P_{i}(t_{j})- \ln P_{i}(t_{1})\right], 
i=1,...,n, \hspace {0.1cm} j=1,...,M.
\]

The CIDRs are observed each minute throughout the trading day.  This corresponds to 
$M=390$ minutes (9:30 am-4:00 pm EST) of trading time for each trading day of the New 
York Stock Exchange, NYSE. In this application 
study, we deal with the CIDR's of two of the most important US market 
indexes: Standard \& Poor's 500 Index (S\&P 500) and the Dow Jones Industrial 
Average (DJ). Also, we consider two individual stocks: General Electric 
Company (GE) and International Business Machines Corporation (IBM). The study period of 
the data consists of three periods in relation to the 2007--2008 financial crisis, denoted as 
 \textit{Before} (06/13/2006-04/10/2007), \textit{During} (11/01/2007-07/28/2008), and 
\textit{After} (01/04/2010-10/1/2010). These periods each contain 270 calendar days.

We investigate the performance of the 
market indexes, S\&P 500 and DJ, in predicting GE and IBM for the three periods.  Understanding such relationships is imperative for developing financial portfolios as many strategies consist of balancing buying/shorting certain stocks with buying/shorting market indices \citep{nicholas2000market}.  
We fit four different models; the first two are our discussed models, AFFR, one based on using the full $X_i(t)$ to predict $Y_i(t)$ (AFFR) and one where only the current and past values are used (AFFR Pre) as described in Section \ref{s:pre}. The other two methods are the linear models.  The first linear model uses an FPCA on both the outcome and predictor (5 PCs for both) and then fits a multivariate linear model, while the second linear model only uses FPCA on the predictor (5 PCs) \citep{KRBook}. To evaluate the prediction performance for each period we split each period into 3 equal folds and use a type $K$-fold cross-validation.  The model is fit on two folds, while prediction is then evaluated on the third.  We use the Gaussian kernel from Section \ref{s:sims} and the smoothing parameter selected via Generalized Cross-Validation.  Prediction performance is then averaged over the 3 folds.  To provide a more readily interpretable metric for prediction performance, we use the same RPE metric given in Section \ref{s:sims}, which denotes the relative performance of a model with respect to a mean only model.  A value of $1$ means perfect prediction, while a value of $0$ indicates that the model is doing no better than just using the mean.  The results are summarized in Table \ref{t:app}.

\begin{figure}
\centering
\includegraphics[width=.99\textwidth]{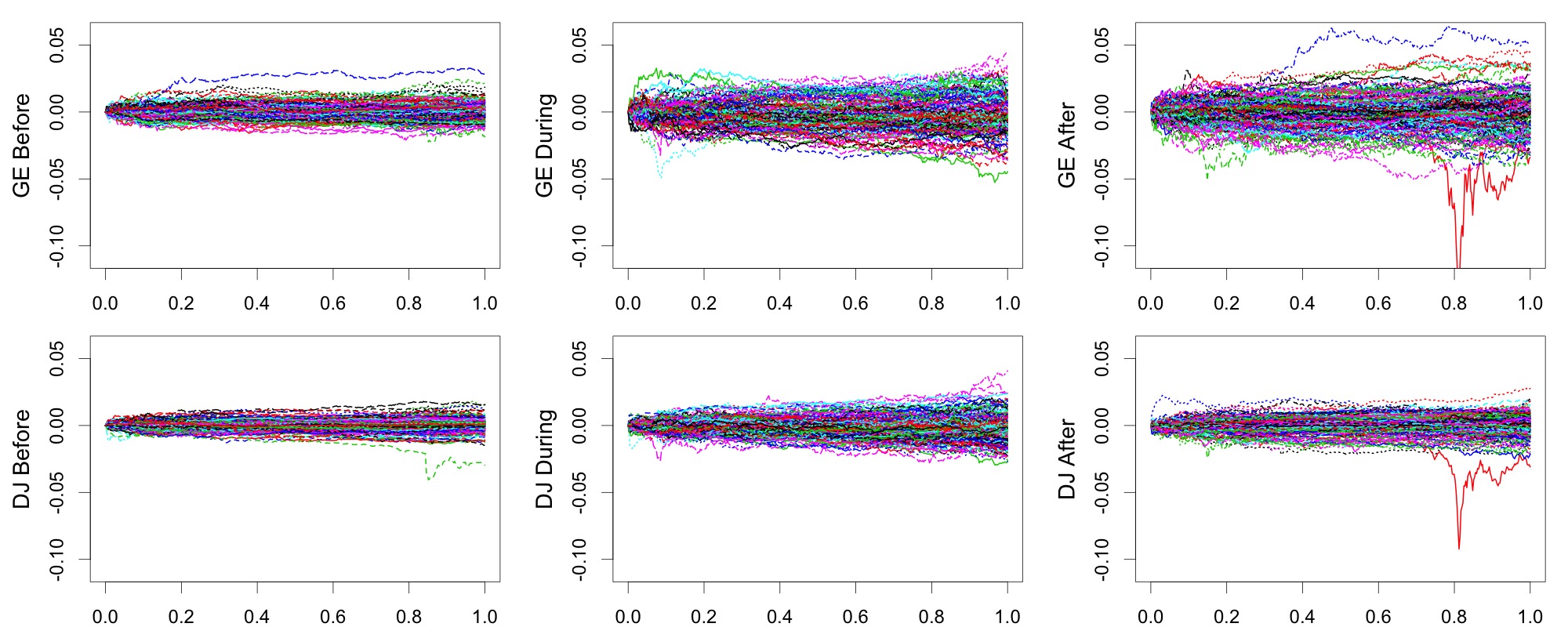}
\caption{Plots of the intraday cumulative returns for the Dow Jones Index (top) and General Electric (bottom) \textit{before} (left), \textit{during} (middle), and \textit{after} (right) the 2008 financial crisis.\label{f:app}}
\end{figure}

As we can see, all models perform better during and after the crisis.  This suggests that the behavior of the market had not returned to its pre-crisis characteristics.  Looking at Figure \ref{f:app}, we can clearly see that the volatility increases during and after the financial crises.  This suggests that the overall ``market" effect on the stocks is stronger during periods of high-volatility.  When comparing the four different models, the linear models do nearly the same, which is to be expected since 5 PCs explains over $90\%$ of the variability of the stocks.  The AFFR model is not too far behind, but does noticeably worse in every setting.  This suggests that the relationship between the discussed stocks and the indices is approximately linear; if there are any nonlinear relationships then they are either very minor deviations from linearity or are not well captured by an additive structure.  The results of AFFR using only current and past values of $X_i(t)$ to predict $Y_i(t)$ (AFFR Pre) does substantially worse before the crises.  Interestingly, during and after its performance is closer to AFFR, though some relationships it still does not capture well. Thus suggests that, unsurprisingly, knowing the future values of $X_i(t)$ is very helpful for predicting currentvalues of $Y_i(t)$, though this is obviously impractical.  During the financial crises, many stocks are likely being driven by large market level effects.  In this setting, AFFR Pre, does quite well, even beating AFFR slightly in some settings, suggesting that the simpler structure has actually helped with prediction.

\begin{table}[ht]
\centering
\resizebox{\columnwidth}{!}{
\begin{tabular}{|c|cccc|cccc|cccc|} \hline
Period & \multicolumn{4}{c|}{Before}&  \multicolumn{4}{|c|}{During}&  \multicolumn{4}{|c|}{After}\\
  \hline
 Model & AFFR & AFFR Pre & LM Red & LM Full & AFFR & AFFR Pre & LM Red & LM Full & AFFR & AFFR Pre & LM Red & LM Full \\ 
 \hline
 GE on DJ & 0.133 & 5.124e-06 & 0.191 & 0.191 & 0.459 & 0.311 & 0.536 & 0.548 & 0.500 & 0.421 & 0.501 & 0.512 \\ 
 GE on SP & 1.325e-07 & 4.216e-14 & 0.184 & 0.183 & 0.273 & 0.253 & 0.458 & 0.472 & 0.510 & 0.436 & 0.487 & 0.497 \\ 
 IBM on DJ & 0.092 & 1.645e-03 & 0.182 & 0.184 & 0.274 & 0.350 & 0.486 & 0.495 & 0.364 & 0.011 & 0.402 & 0.412 \\ 
 IBM on SP & 0.079 & 1.251e-11 & 0.180 & 0.180 & 0.213 & 0.272 & 0.373 & 0.384 & 0.296 & 0.009 & 0.343 & 0.351 \\   
   \hline
\end{tabular}
}
\caption{Prediction performance of four models: AFFR (our model), AFFR Pre (modifies domain to avoid using future values), LM Red (linear model with PCA in both the outcome and predictor), and LM Full (linear model with PCA on the predictor only).  The top row corresponds to predicting GE based on DJ, the second corresponds to prediction GE from SP, and so on.  Each number denotes the relative increase in out-of-sample prediction performance over a mean only model, with $100\%$ denoting perfect prediction and $0\%$ denoting no increase over just using the mean.  \label{t:app}}
\end{table}

\section{Conclusions and Future Work}
In this paper we have presented a new RKHS framework for estimating an additive function-on-function regression model, that is better able to account for complex nonlinear dynamics in functional regression models than classic linear models.  We showed that the estimator is minimax in the sense that it achieves an optimal rate of convergence in terms of prediction error.  In addition, computing the estimate is computationally efficient, especially if dimension reduction is incorporated.

Nonlinear models for functional data have recently received a great deal of attention, however, there are still a number of interesting questions that remain open.  One that is especially relevant to the work presented here concerns further statistical properties of the estimate, $\hat g$.  In particular, convergence rates of $\hat g$ as well as its asymptotic distribution would be especially interesting for quantifying the estimation uncertainty in practice.  Using such tools, one could also construct confidence/prediction bands, which would be of great use in practice.

Another nontrivial extension would be to curves that are observed sparsely.  Nonlinear models in FDA often require that the curves be observed or at least consistently estimated.  However, for some data this is unrealistic and there is a great deal of uncertainty related to imputing the curves. 

Lastly, extensions to more complex settings would also be of interest. For example, the handling of more complex domains, e.g. space or space-time.  In these cases, the minimax rates usually depend on the dimension of the domain.  Another important extension would be to functional binary or categorical outcomes (as opposed to quantitative) would be of interest as one must incorporate tools from functional glms.    

\appendix

\section{Proof of Theorem \ref{t:main}}
\subsection{Excess Risk }
We begin by expressing the excess risk in an alternative form. Recall that 
$k(t,s,x; t',s',x')$ is the kernel function used to define the RKHS, $\mbK$. 
This kernel can be viewed as the kernel of an integral operator, $\mcL_k$, which 
maps $L^2( [0,1]^2 \times \mbR) \to \mbK \subset L^2( [0,1]^2 \times \mbR) $. 
In particular
\[
(\mcL_k f)(t,s,x) = \int \int \int k(t,s,x; t', s', x') f(t',s',x') \ dt' ds' 
dx'.
\]
 From here on, for simplicity, we will denote $L^2( [0,1]^2 \times \mbR)$ as 
simply $L^2$. By Assumption \ref{a:main}, $\mathcal{L}_{k}$ is a  
positive definite, compact operator, which is also self-adjoint in the sense that $\langle f, \mcL_k g \rangle =  \langle \mcL_k f,  g \rangle$, for any $f$ and $g$ in $L^2$.  
We can therefore define a  
square-root of $\mcL_k$, denoted as $\mcL_k^{1/2}$ that satisfies 
\[
 f_{1} \in \mbK \Longleftrightarrow \mathcal{L}_{k}^{-\frac{1}{2}} f_{1} \in 
L^{2} \hspace {0.3cm}\text{and} \hspace {0.3cm} \mathcal{L}_{k}^{\frac{1}{2}} 
f_{2}\in \mbK \Longleftrightarrow f_{2} \in L^{2}.
\]
Note that if $\mcL_k$ has a nontrivial 
null space, then $\mcL_k^{-1/2}$ can 
still be well defined since assuming $f \in \mbK$ means that $f$ is orthogonal 
to the null space of $\mcL$. Recall that one can also move between the $\mbK$ 
and $L^2$ inner product as follows
\[
\langle f, g \rangle_\mbK = \langle 
\mcL_k^{-1/2} f, \mcL_k^{-1/2} g \rangle_{L^2}=\langle f, \mcL_k^{-1} g \rangle_{L^2}.
\]
We refer the interested reader to \citet{kennedy2013hilbert} for more 
details.

Let $\hat g$ denote our estimate of the true function, $g$. We then define the 
following 
\[
 k^{\frac{1}{2}}_{t,s, X(s)} = \mathcal{L}_{k}^{-\frac{1}{2}} k_{t,s, 
X(s)},\quad 
 h= \mathcal{L}_{k}^{-\frac{1}{2}} g \hspace {0.3cm} \text{and} \hspace 
{0.3cm} \hat h = \mathcal{L}_{k}^{-\frac{1}{2}} \hat g.  
\]
Using the reproducing property, we have that 
\[ 
 g(t,s,X(s)) 
 = \langle k^{\frac{1}{2}}_{t,s, X(s)}, h \rangle_{L^{2}}  
\quad \text{and} \quad
\hat g(t,s,X(s)) 
= \langle k^{\frac{1}{2}}_{t,s, X(s)}, \hat h \rangle_{L^{2}}.  
\]
Now define the random operator, $T: L^2\rightarrow L^2$ as 
\[ 
T_{t,s,s'} = k^{1/2}_{t,s,X_{n+1}(s)} \otimes k^{1/2}_{t,s',X_{n+1}(s')},
\]
where $\otimes$ denotes the tensor product, and the resulting object is 
interpreted as an operator:
\[
T_{t,s,s'}(f) = k^{1/2}_{t,s,X_{n+1}(s)} \langle k^{1/2}_{t,s',X_{n+1}(s')}, f 
\rangle_{L^2}.
\]
We also define a second operator, which integrates out $t$, $s$, and $s'$, and 
takes an expectation over $X_{n+1}$:
\begin{align}
 C = \E\left[ \int_0^{ 1}\int_0^{ 1} \int_0^{ 1} T_{t,s,s'} ds ds' dt 
\right]. \label{e:C}  
\end{align}
Note that $C$ is a symmetric, positive definite, compact operator, and thus has 
a spectral decomposition
\begin{align}
C = \sum_{k=1}^{\infty}\rho_{k} (\phi_{k} \otimes \phi_{k}), 
\label{e:C_eigen}  
\end{align}
where $\rho_{k} \geq 0$ and $ \phi_{k} \in L^2$ are, respectively, the 
eigenvalues and eigenfunctions of $C$. This decomposition will be used later 
on.
\\
\\
As we said before, denote by $ {\E^*}$ the expected value conditioned on the 
data $(X_{1},Y_{1}),...,(X_{n},Y_{n})$.
The excess risk can be written as
\begin{align*} 
 \Re_{n} 
&= {\E^*}\int_0^{ 1} \left(\int_0^{ 1}\left[\hat g(t,s,X_{n+1}(s))- g(t,s,X_{n+1}(s))) 
\right] ds \right)^{2} dt\\
  &= {\E^*} \int_0^{ 1} \left(\int_0^{ 1}\left[ \langle 
k^{\frac{1}{2}}_{t,s, X_{n+1}(s)}, \hat h_{ \lambda} \rangle_{L^{2}} - \langle 
k^{\frac{1}{2}}_{t,s, X_{n+1}(s)}, h \rangle_{L^{2}} \right]ds \right)^{2}dt \\
  &= {\E^*} \int_0^{ 1} \left(\int_0^{ 1}\langle 
k^{\frac{1}{2}}_{t,s, X_{n+1}(s)}, \hat h_{ \lambda}-h \rangle_{L^{2}} ds 
\right)^{2} dt\\
 &= {\E^*} \int_0^{ 1}\int_0^{ 1} \int_0^1 \langle k^{\frac{1}{2}}_{t,s, 
X_{n+1}(s)}, \hat h-h \rangle_{L^{2}} 
 \langle k^{\frac{1}{2}}_{t,s', X_{n+1}(s')}, \hat h-h \rangle_{L^{2}}
 ds ds' dt\\
 & = {\E^*} \int_0^{ 1}\int_0^{ 1} \int_0^1 \langle T_{t,s,s'} (\hat h - h), 
\hat h-h \rangle_{L^{2}} ds ds' dt \\
  &= \langle C(\hat h - h),\hat h - h \rangle_{L^2} = \|\hat 
h-h\|^{2}_{C}. 
\end{align*} 
Thus, the excess risk can be expressed as sort of a weighted $L^2$ norm, where 
the operator $C$ defines the weights, which is composed of the kernel and the 
distribution of $X_{n+1}$.

\subsection{Re-expressing the Estimator }
In this section we define an alternative form for the estimator $\hat g$, which 
was given in Section \ref{s:comp}. In particular, instead of using the 
reproducing property, we will write down the estimator using operators. To do 
this, we will take derivatives of $RSS_\lambda(g)$ with respect to $g$. Since 
these are functions, we mean the Fr\'echet derivative or strong derivative. 
Note that $RSS_\lambda(g)$ is a convex differentiable functional over $\mbK$. 
However, so that we are working with $L^2$ instead of $\mbK$, we use $\widetilde{RSS}_\lambda(h) : = RSS_\lambda(\mcL^{1/2}_k h)$, where $h = \mcL^{-1/2}_k g$:  
\[
\widetilde{RSS}_\lambda(h) = 
\sum_{i=1}^n \int_0^1 \left(
Y_i(t) - \int_0^1 \langle h, k^{1/2}_{t,s,X_i(s)} \rangle_{L^2} ds 
\right)^2 dt + 
\lambda \|h \|^2_{L^2}.
\]
Now 
$\widetilde{RSS}_\lambda(h)$ is a convex differentiable functional over $L^2$. 
Thus, when taking the derivative, we are using the topology of  $L^2$ not 
 $\mbK$. 

To take the derivative of $\widetilde{RSS}_\lambda(h)$ we first focus on the penalty, 
which is easier. We have that
\[
\frac{\partial}{\partial h} \| h\|^2_{L^2} = 2 h. 
\]
Turning to the first term in $\widetilde{RSS}_\lambda(h)$ we first define the 
empirical quantities
\[
T_{i; t,s,s'} = k^{1/2}_{t,s,X_{i}(s)} \otimes k^{1/2}_{t,s',X_{i}(s')}
\]
and
\begin{align}
C_n = \frac{1}{n} \sum_{i=1}^n \int_0^{ 1}\int_0^{ 1} \int_0^{ 1} T_{i;t,s,s'} 
ds ds' dt . \label{e:cn}
\end{align}
Now we can apply a chain rule to obtain
\begin{align*}
 & \frac{\partial }{\partial h} \left[ \frac{1}{n}\sum_{i=1}^{n} \int_0^{ 
1}\left (Y_{i}(t)- \int_0^1 \langle h, k^{1/2}_{t,s,X_i(s)} \rangle_{L^2} ds  \right)^2 dt\right] \\
 &= -\frac{2}{n}\sum_{i=1}^{n} \int_0^{ 1} \int_0^{ 1} Y_{i}(t) 
k^{\frac{1}{2}}_{t,s, X_{i}(s)}ds dt
 + 2 C_n h. 
\end{align*}
For notational simplicity, define
\[
\Gamma_{k^{1/2},Y} = 
\frac{1}{n}\sum_{i=1}^{n} \int_0^{ 1} \int_0^{ 1} Y_{i}(t) 
k^{\frac{1}{2}}_{t,s, X_{i}(s)}ds dt.
\]
So, we finally have that
\[
\frac{\partial}{\partial h} \widetilde{RSS}_\lambda(h) 
= - 2\Gamma_{k^{1/2},Y} + 2 C_n h + 2 \lambda h ,
\]
which yields the estimate
\begin{align}
\hat h = ( C_n + \lambda I)^{-1} \Gamma_{k^{1/2},Y}, \label{e:hhat}
\end{align}
where $I$ is the identity operator.

\subsection{Proof of Theorem \ref{t:main} - Controlling Bias}\label{s:bias}
 Using Assumption \ref{a:main} we can express
 \[
 Y_i(t) = \int \langle k^{1/2}_{t,s,X_i(s)}, h \rangle_{L^2} + \vep_i(t). 
 \]
and we therefore have that 
\[
\Gamma_{k^{1/2},Y} = C_n(h) + f_n
\]
where 
\[
 f_{n}=\frac{1}{n} \sum_{i=1}^{n} \int_0^{ 1} \int_0^{ 1} \epsilon_{i}(t) 
k^{\frac{1}{2}}_{t,s, X_{i}(s)} ds dt.
\]
This implies that $\hat h$ from \eqref{e:hhat} can be expressed as
\[
\hat h = (C_n + \lambda I)^{-1}C_n(h) + (C_n + \lambda I)^{-1}f_n.
\]
We introduce an intermediate quantity, $h_\lambda$, which is given by
\[
h_\lambda = (C + \lambda I)^{-1} C(h),
\]
where $C$ is defined in \eqref{e:C}. The difference between $h_\lambda$ and $h$ 
represents the bias of the estimator $\hat h$. Balancing this quantity with the 
variance, discussed in the next section, is called the \textit{bias-variance 
trade off} a common term in nonparametric smoothing. Inherently, the idea is 
that to achieve an optimal $\hat h$ we have to balance both the bias and 
variance so that neither one is overly large.

Using the eigenfunctions of $C$ as a basis, we can write
\begin{equation*}
 \begin{split} 
h & = \sum_{k=1}^{\infty} a_{k}\phi_{k}.  
\end{split}
\end{equation*}
Since $C$ and $C+I$ have the same eigenfunctions, it follows that we can 
express
\[
C+\lambda I = \sum_{k=1}^{\infty} (\lambda+ \rho_{k})( \phi_{k} \otimes 
\phi_{k})
\Longrightarrow 
(C+\lambda I)^{-1} = \sum_{k=1}^{\infty} (\lambda+ \rho_{k})^{-1}( \phi_{k} 
\otimes \phi_{k}). 
\]
So we have that $h_\lambda$ can be expressed as
\begin{equation*}
 \begin{split} 
h_{\lambda}=(C+\lambda I)^{-1} C(h) 
= \sum_{k=1}^{\infty} \frac{ a_{k} \rho_{k}}{\lambda + \rho_{k}}\phi_{k}. 
\end{split}
\end{equation*}
So the difference, $h_\lambda - h$ can be written as
\begin{align}
 h_{\lambda}-h= - \sum_{k=1}^{\infty} \frac{ \lambda a_{k} 
}{\lambda + \rho_{k}}\phi_{k}. \label{e:h_dif}
\end{align}
The bias is therefore given by
\[ 
 \| h_{\lambda}-h \|^{2}_{C} = \sum_{k=1}^{\infty} \frac{ \lambda^{2} 
a_{k}^{2}\rho_{k} }{(\lambda + \rho_{k})^{2}}
\leq \lambda^2 \max_{k\geq 1} \frac{\rho_{k}}{(\lambda+\rho_{k})^2} 
\sum_{k=1}^{\infty} a_{k}^2
= \lambda^2 \| h\|^2_{L^2} \max_{k\geq 1} 
\frac{\rho_{k}}{(\lambda+\rho_{k})^2}.
\]
It is easy to verify that the maximum of $F(x) = x/(\lambda+x)^2$ is achieved at $x = \lambda$ with the maximum value being $\frac{1}{4\lambda}$. We can therefore bound the bias 
as
\[
\| h_\lambda - h \|^2_C
\leq \frac{\lambda \| h\|^2_{L^2}}{4}. 
\]
In the statement of Theorem \ref{t:main} we assume that $\lambda \asymp n^{2r/(2r+1)}$,
which implies that the bias is of the order 
$n^{-\frac{2r}{2r+1}} O(1)$. We will show in the next section that the variance 
of our estimate achieves the same order.

\subsection{Proof of Theorem \ref{t:main} - Controlling Variability} 
\label{s:var}
Controlling the variability of the estimates, $\| \hat h - h_\lambda\|_C$ follows 
similar arguments as controlling the bias. However, there are many more terms 
which must be analyzed separately. In particular, we decompose $\hat h - 
h_\lambda$ into five separate components:
\begin{align}
\hat h - h_\lambda = T_1 + T_2 + T_3 + T_4 + T_5, \label{e:T}
\end{align}
where the $T_i$ terms are given by
\begin{align*}
 T_{1}& = (C+\lambda I)^{-1} C (h_{\lambda}-h),\\
 T_{2}&= \lambda (C+\lambda I)^{-2} C(h),\\
 T_{3}&= -(C+\lambda I)^{-1} f_{n},\\
 T_{4}&= (C+\lambda I)^{-1} (C_{n}-C) (h_{\lambda}-h),\\
 T_{5}&= (C+\lambda I)^{-1} (C-C_{n})( h_{\lambda}-\hat h). 
\end{align*}
While a bit tedious, it only requires linear algebra and repeated calls to the 
definitions of $\hat h$ and $h_\lambda$ to verify \eqref{e:T}, we thus omit the 
details here. We now develop bounds for each term, $\|T_i\|_C$, separately. For the first 
four, it turns out to be convenient to bound $\| C^\nu T_i\|_{L^2}$ for $0 < \nu 
\leq 1/2$, as these bounds will be needed for the final term $T_5$. Notice that  when $\nu =1/2$ we have $\| C^\nu T_i\|_{L^2} = \| T_i \|_C$.

\begin{enumerate}
\item Using the eigenfunctions of $C$ to express $h_\lambda - h$ as in 
\eqref{e:h_dif}, we get that
\[
T_1 = - \sum_{k=1}^{\infty} \frac{ \lambda a_{k} \rho_k}{(\lambda + 
\rho_{k})^2}\phi_{k}.
\]
We then have that
\[
\| C^\nu T_1 \|^2_{L^2} = \sum_{k=1}^{\infty} \frac{ \lambda^{2} 
a_{k}^{2}\rho_{k}^{2(1+\nu)} }{(\lambda + \rho_{k})^{4}}
 \leq \lambda^2 \max_{k\geq 1} 
\frac{\rho_{k}^{2(1+\nu)}}{(\lambda+\rho_{k})^4} \| h \|^{2}_{L^2}.
\]
Again, it is a basic calculus exercise to show that 
\[
\max_{k\geq 1} \frac{\rho_{k}^{2(1+\nu)}}{(\lambda+\rho_{k})^4}
\leq \frac{\left(\lambda 
\frac{1+\nu}{1-\nu}\right)^{2(1+\nu)}}{\left(\lambda+\lambda 
\frac{1+\nu}{1-\nu}\right)^4}
= \frac{(1-\nu)^{2(1-\nu)} (1+\nu)^{2(1+\nu)}}{16} \frac{1}{\lambda^{2-2\nu}}.
\]
We thus have the bound
\begin{align}
\| C^\nu T_1\|^2_{L^2} \leq c \lambda^{2\nu} \|h\|^2_{L^2}, \label{e:T1}
\end{align}
where $c$ is a constant that depends only on $\nu$.

\item Using the same arguments as in the previous step, we have that
\begin{align}
 \| C^\nu T_{2} \|^{2}_{L^2} 
 = \sum_{k=1}^{\infty} \frac{ \lambda^{2} a_{k}^{2}\rho_{k}^{2(1+\nu)} 
}{(\lambda + \rho_{k})^{4}}
\leq \lambda^2 \max_{k\geq 1} \frac{\rho_{k}^{2(1+\nu)}}{(\lambda+\rho_{k})^4} 
\sum_{k=1}^{\infty} a_{k}^2
\leq c \lambda^{2\nu} \|h\|^2_{L^2}. \label{e:T2}
\end{align}

\item Turning to $T_3$, we apply Lemma \ref{lemma1} with $0 < \nu \leq 1/2$ to 
obtain
\[
\| C^\nu T_3\|^2_{L^2} = \| C^{\nu}(C+\lambda I)^{-1} f_{n} \|_{L^2} = 
\frac{1}{n \lambda^{1-2 \nu+ 1/2r} } O_p(1),
\]
where $r$ is defined as in Assumption \ref{a:main}. By the statement of Theorem \ref{t:main} 
it follows that $n \lambda ^{1+\frac{1}{2r}}$ tends to a nonzero constant, 
meaning that
\begin{align}
\frac{1}{n \lambda^{1-2\nu + \frac{1}{2r}}}
\asymp \lambda^{2 \nu} \to 0, \label{e:nu_rate}
\end{align}
since $\lambda \to 0$. Thus we have that $\| C^\nu T_3\|_{L^2} = 
O_p(\lambda^\nu)$.


\item To bound $T_4$ we first fix a second value $\nu > \nu_2 > 0$ that satisfies 
$2r(1-2\nu_2) > 1$, or equivalently $\nu_2 < (2r -1)/4r$, as well as $4r(2\nu_2 + 2 \nu) > 1$, which is possible as 
long as $r > 1/2$ (Assumption \ref{a:main}).  We now apply a basic operator inequality
\begin{equation*}
 \begin{split} 
 \| C^\nu T_{4}\|_{L^2}&= \| C^{\nu}(C+\lambda I)^{-1} (C_{n}-C) 
(h_{\lambda}-h)\|_{L^2}\\
 & \leq \| C^{\nu}(C+\lambda I)^{-1} (C_{n}-C) C^{-\nu_2}\|_{op}\| 
C^{\nu_2}(h_{\lambda}-h) \|_{L^2}.
\end{split}
\end{equation*}
and then apply Lemmas \ref{lemma3} and \ref{lemma4} to obtain
\begin{align} 
 \| C^\nu T_{4}\|_{L^2}^2& \leq O_{p}\left(\left( n \lambda^{1-2\nu+\frac{1}{2r}} 
\right)^{-1} \right) O_{p}\left( \lambda^{2\nu_2}\right)
 = \frac{1}{n \lambda^{1-2\nu+\frac{1}{2r}}}o_{p}\left(1\right),
\end{align}
since $\nu_2 > 0$ and $\lambda \to 0$. Using \eqref{e:nu_rate} we conclude that 
$\|C^\nu T_4\|_{L^2} = O_p(\lambda^\nu)$.

\item The last term is the most involved to bound and the reason why the 
previous four bounds involved $C^\nu T_i$. We begin by expressing
\begin{align*}
\| T_{5}\|_{C}&= \| C^{\frac{1}{2}}(C+\lambda I)^{-1} (C-C_{n})( 
h_{\lambda}-\hat h)\|_{L^2}\\
& \leq \| C^{\frac{1}{2}}(C+\lambda I)^{-1} (C-C_{n})C^{-\nu} \|_{op}\| 
C^{\nu}( h_{\lambda}-\hat h) \|_{L^2}.
\end{align*} 
Here $\nu>0$ is chosen to satisfy $2r(1-2\nu) > 1$. Applying Lemma \ref{lemma3}, we have 
that
\[
\|T_5\|_C^2 \leq \frac{1}{n \lambda^{1/2r}}O_p(1) \| C^\nu(h - \hat 
h)\|_{L^2}^2.
\]
We have now, in some sense, looped back and are dealing with the term $h - \hat 
h$. Using \eqref{e:T} we have
\begin{align}
\| C^\nu( h - \hat h)\|_{L^2} \leq \| C^\nu T_1\|_{L^2} + 
\| C^\nu T_2\|_{L^2}+
\| C^\nu T_3\|_{L^2} +
\| C^\nu T_4\|_{L^2} +
\| C^\nu T_5\|_{L^2}. \label{e:Tnu}
\end{align}
The first four terms we already have bounds for, so we need only focus on the 
last, which again, has looped back to our original term. We now apply Lemma 
\ref{lemma2} to obtain
\[
\|C^\nu T_5\|_{L^2}^2
= \| C^\nu (C+\lambda I)^{-1} (C- C_n)(h - \hat h)\|_{L^2}^2
\leq \frac{1}{n \lambda^{1-2\nu + \frac{1}{2r}}} \|C^\nu (h - \hat h)\|_{L^2}^2.
\]
Combining the above with \eqref{e:Tnu} we have that
\[
\| C^\nu( h - \hat h)\|_{L^2} \left( 1 - \frac{1}{n \lambda^{1-2\nu + 
\frac{1}{2r}}} \right) \leq \| C^\nu T_1\|_{L^2} + 
\| C^\nu T_2\|_{L^2}+
\| C^\nu T_3\|_{L^2} +
\| C^\nu T_4\|_{L^2}.
\]
Using \eqref{e:nu_rate} it thus follows that
\[
\| C^\nu( h - \hat h)\|_{L^2} = O_p(1) ( \| C^\nu T_1\|_{L^2} + 
\| C^\nu T_2\|_{L^2}+
\| C^\nu T_3\|_{L^2} +
\| C^\nu T_4\|_{L^2} ),
\]
and applying steps 1-4 we get that
\[
\| C^\nu( h - \hat h)\|_{L^2} = O_p(\lambda^\nu) =o_p(1),
\]
and we finally have that
\[
\|T_5\|_{C} = \frac{1}{n \lambda^{\frac{1}{2r}}} o_p(1).
\]
\end{enumerate}
We can now combine Steps 1-4, taking $\nu = 1/2$, with step 5 to finally 
conclude that
\[
\Re_n^2 = \| \hat h - h \|_C^2 = \lambda^{2} O_p(1)
= n^{-\frac{2r}{2r+1}} O_p(1).
\]
Combined with the results of Section \ref{s:bias}, this concludes the proof.

\subsection{Auxiliary Lemmas}

Here we state four lemmas which are generalizations of ones used in 
\citet{cai2012minimax} and \citet{wang2015optimal}.

\begin{lemma} \label{lemma1}
	If Assumption \ref{a:main} holds then for any $0 \leq \nu \leq 
\frac{1}{2}$
	\begin{equation*}
	\begin{split} 
	\| C^{\nu}(C+\lambda I)^{-1} f_{n} \|_{L^2} = O_{p}\left(\left( n 
\lambda^{1-2\nu+\frac{1}{2r}} \right)^{-\frac{1}{2}}\right).
	\end{split}
	\end{equation*}
\end{lemma}

\begin{lemma} \label{lemma2}
	Let Assumption \ref{a:main} hold. 
	Then for any $\nu >0$ such that $2r(1-2\nu) >1$, we have that
	\begin{equation*}
	\begin{split} 
	\| C^{\nu}(C+\lambda I)^{-1} (C_{n}-C) C^{-\nu} \|_{op} 
=O_{p}\left(\left( n \lambda^{1-2\nu+\frac{1}{2r}} 
\right)^{-\frac{1}{2}}\right),
	\end{split}
	\end{equation*}
	where $\| . \|_{op}$ represents the usual operator norm i.e., $\| 
A\|_{op} =\sup_{h: \| h \|_{L^2}=1} \| Ah \| $.
\end{lemma}

\begin{lemma} \label{lemma3}
Let Assumption \ref{a:main} hold and fix $0 < \nu< \nu_2$ to be any two values that satisfy $2 r(1-2 
\nu) > 1$ and $4r(\nu_2 + \nu) > 1$, then we have that
\[
	\| C^{\nu_2}(C+\lambda I)^{-1} (C_{n}-C) C^{-\nu} \|_{op} 
=O_{p}\left(\left( n \lambda^{ 1 - 2 \nu_2+ \frac{1}{2r}} 
\right)^{-\frac{1}{2}}\right).
\]
\end{lemma}

\begin{lemma} (\citealp[Lemma 1]{cai2012minimax}) \label{lemma4}
	For any $ 0<\nu<1$,
	\begin{equation*}
	\begin{split} 
	\| C^{\nu}(h_{\lambda}-h) \|_{L^2}\leq 
(1-\nu)^{1-\nu}\nu^{\nu}\lambda^{\nu}\| h \|_{L^2}.
	\end{split}
	\end{equation*}
\end{lemma}

 \begin{lemma}\label{lemma5} 
Fix $\nu>0$ and $\nu_2>0$ such that $4r(\nu_2 + \nu) > 1$.  If there exist constants $0<c_{1}< c_{2}< \infty$ such that $ c_{1} 
k^{-2r}<s_{k}<c_{2} k^{-2r}$, then there exist constants $c_{3},c_{4}>0$ 
depending only on $c_{1}, c_{2}$ such that 
\begin{equation*}
 \begin{split} 
c_{4} \lambda^{\frac{-1}{2r} - 1 + 2 \nu_2}\leq 
\sum_{j=1}^{\infty}\frac{s_{j}^{2\nu_2+2\nu}}{(\lambda+s_{j})^{1+2\nu}} \leq 
c_{3}(1+\lambda^{\frac{-1}{2r} - 1 + 2 \nu_2}).
\end{split} 
\end{equation*}
\end{lemma}

\subsection*{Proof of Lemma \ref{lemma1}}
Recall that 
\[
 f_{n}=\frac{1}{n} \sum_{i=1}^{n} \int_0^{ 1} \int_0^{ 1} \epsilon_{i}(t) 
k^{\frac{1}{2}}_{t,s, X_{i}(s)} ds dt.
\]
Using Parseval's identity we have that
\begin{align*}
& \| C^{\nu}(C+\lambda I)^{-1}f_{n} \|^{2}_{L^2} 
 = \sum_{k=1}^\infty \frac{\rho_k^{2\nu}}{(\lambda+ \rho_k)^2}\langle f_n, 
\phi_k \rangle^2.
\end{align*}
Taking expected values yields
\begin{align*}
 & \E \| C^{\nu}(C+\lambda I)^{-1}f_{n} \|^{2}_{L^2}
 = \frac{1}{n} 
\sum_{k=1}^{\infty}\frac{\rho_{k}^{2\nu}}{(\lambda+\rho_{k})^{2}} \E \left ( 
\int_0^{ 1} \int_0^{ 1} \epsilon(t) \langle k^{\frac{1}{2}}_{t,s, X(s)} , 
\phi_{k} \rangle_{L^2} ds dt \right)^{2}.
\end{align*}
By Jensen's inequality we have
\normalsize
\begin{align*}
  \left ( \int_0^{ 1} \int_0^{ 1} \epsilon(t) \langle k^{\frac{1}{2}}_{t,s, 
X(s)} , \phi_{k} \rangle_{L^2} ds dt \right)^{2} 
 & \leq \int_0^{ 1} \left( \int_0^{ 1} \epsilon(t) \langle 
k^{\frac{1}{2}}_{t,s, X(s)} , \phi_{k} \rangle_{L^2} ds\right)^{2}dt  \\
&=\int_0^{ 1} \int_0^{ 1} \int_0^{ 1} \epsilon^{2}(t) \langle 
k^{\frac{1}{2}}_{t,s, X(s)} , \phi_{k} \rangle_{L^2} \langle 
k^{\frac{1}{2}}_{t,s^{*}, X(s^{*})} , \phi_{k} \rangle_{L^2}ds ds^{*} dt.
\end{align*}
Using the assumed independence between $\epsilon$ and $X$, as well as the 
assumption that $ \E (\epsilon^{2}(t)) \leq M$, where $M$ is a constant, we 
obtain 
\normalsize
\begin{align*}
 \E \left ( \int_0^{ 1} \int_0^{ 1} \epsilon(t) \langle 
k^{\frac{1}{2}}_{t,s, X(s)} , \phi_{k} \rangle_{L^2} ds dt \right)^{2}
 & \leq M \E \left( \int_0^{ 1} \int_0^{ 1} \int_0^{ 1} \langle 
k^{\frac{1}{2}}_{t,s, X(s)} , \phi_{k} \rangle_{L^2} \langle 
k^{\frac{1}{2}}_{t,s^{*}, X(s^{*})} , \phi_{k} \rangle_{L^2}ds ds^{*} 
dt\right)\\
     & = M 
\langle C(\phi_{k}) ,\phi_{k}\rangle 
     =M 
\rho_{k}. 
\end{align*}
Since $0\leq \nu \leq \frac{1}{2}$ and both $\rho_k$ and $\lambda$ are positive, 
we can obtain the bound
\[ 
 \E \| C^{\nu}(C+\lambda I)^{-1}f_{n} \|^{2}_{L^2}\leq 
\frac{M}{n}\sum_{k=1}^{\infty}\frac{\rho_{k}^{2\nu + 1}}{(\lambda+\rho_{k})^{2}}
 \leq \frac{M}{n \lambda^{1-2\nu}}\sum_{k=1}^{\infty}\frac{\rho_{k}^{2\nu + 
1}}{(\lambda+\rho_{k})^{1+2\nu}}.
    \]
Now we apply Lemma \ref{lemma5} with $\nu_2 = 1/2$  to obtain
\begin{equation*}
 \begin{split} 
 \E \| C^{\nu}(C+\lambda I)^{-1}f_{n} \|^{2}_{L^2} & \leq \frac{c^{*}}{n \lambda 
^{1-2\nu+\frac{1}{2r}}},
 \end{split}   
\end{equation*}
where $c^{*}$ is a constant. An application of Markov's inequality completes 
the proof.

\subsection*{Proof of Lemma \ref{lemma2}}
By definition
\normalsize
\begin{equation*}
 \begin{split} 
 \| C^{\nu}(C+\lambda I)^{-1} (C_{n}-C) C^{-\nu} \|^{2}_{op} & = \sup_{f: 
\| f \|_{L^2}=1} \| C^{\nu}(C+\lambda I)^{-1} (C_{n}-C) C^{-\nu} f 
\|^{2}_{L^2}.
\end{split} 
\end{equation*}
Fix $f \in L^2$ such that $\| f\|_{L^2} = 1$.  We can expand $f$ as 
\begin{equation*}
 \begin{split} 
 f& =\sum_{k=1}^{\infty}f_{k} \phi_{k},
\end{split} 
\end{equation*}
%
%
By Parseval's identity we have
\begin{equation*}
 \begin{split} \| C^{\nu}(C+\lambda I)^{-1} (C_{n}-C) C^{-\nu}f\|^{2}& = 
\sum_{j=1}^{\infty}\left[\frac{\rho_{j}^{\nu}}{\rho_{j}+\lambda}\sum_{k=1}^{\infty}f_{k}\rho_{k}^{-\nu} \langle 
 (C_{n}-C)\phi_{k},\phi_{j} 
\rangle_{L^2}\right]^{2}.
\end{split} 
\end{equation*}
Applying the Cauchy-Schwartz inequality and using the fact that $\|f\|_{L^2} =1$ we have that 
\begin{align*}
\sum_{k=1}^{\infty}f_{k}\rho_{k}^{-\nu} \langle 
 (C_{n}-C)\phi_{k},\phi_{j} 
\rangle_{L^2}
\leq  \left( \sum_{k=1}^{\infty}\rho_{k}^{-2\nu} \langle 
 (C_{n}-C)\phi_{k},\phi_{j} \rangle_{L^2}^2
 \right)^{1/2}.
\end{align*}
So we can bound the operator norm as
\begin{equation*}
 \begin{split} 
 \| C^{\nu}(C+\lambda I)^{-1} (C_{n}-C) C^{-\nu} \|^{2}_{op} & \leq 
\sum_{k=1}^{\infty} \sum_{j=1}^{\infty}\frac{\rho_{k}^{-2\nu} 
\rho_{j}^{2\nu}}{(\lambda+\rho_{j})^{2}} \langle \phi_{j},(C_{n}-C) \phi_{k} 
\rangle^{2}_{L^2}.
\end{split} 
\end{equation*}
Applying Jensen's equality we get that
\begin{align*}
 & \E \left(\sum_{k=1}^{\infty} \sum_{j=1}^{\infty}\frac{\rho_{k}^{-2\nu} 
\rho_{j}^{2\nu}}{(\lambda+\rho_{j})^{2}} \langle \phi_{j},(C_{n}-C) \phi_{k} 
\rangle^{2}_{L^2}\right)^{\frac{1}{2}} 
 \leq \left(\sum_{k=1}^{\infty} \sum_{j=1}^{\infty}\frac{\rho_{k}^{-2\nu} 
\rho_{j}^{2\nu}}{(\lambda+\rho_{j})^{2}} \E \langle \phi_{j},(C_{n}-C) \phi_{k} 
\rangle^{2}_{L^2}\right)^{\frac{1}{2}}.
\end{align*}
Using the definition of $C_n$ from \eqref{e:cn} we have that
\begin{align*}
  \E\langle \phi_{j}, (C_{n}-C) \phi_{k} \rangle^{2}_{L^{2}} 
 {\leq} \E\langle \phi_{j}, C_{n} \phi_{k} \rangle^{2}_{L^{2}} 
    = \frac{1}{n} \E \left ( \int_0^{ 1} \int_0^{ 
1}\int_0^{ 1} \langle k^{\frac{1}{2}}_{t,s, X(s)} , \phi_{j} \rangle 
\langle k^{\frac{1}{2}}_{t,s^{*}, X(s^{*})} , \phi_{k}\rangle ds ds^{*} dt 
\right)^{2}.
\end{align*}
Note the first inequality follows from the fact that $C$ is the mean of $C_n$ and thus replacing $C$ above with any other quantity cannot decrease it (since it is minimized when using $C$).  One can show this using basic calculus arguments over Hilbert spaces, thus we omit the details here.
By applying Cauchy-Schwartz inequality and Fubini's theorem we have
\begin{align*}
 &\E\langle \phi_{j}, (C_{n}-C) \phi_{k} \rangle^{2}_{L^{2}}\\
 & \leq \frac{1}{n} \E \left( \left[\int_0^{ 1} \left( \int_0^{ 1} 
\langle k^{\frac{1}{2}}_{t,s, X(s)} , \phi_{j} \rangle ds 
\right)^{2}dt\right] \left[\int_0^{ 1} \left(\int_0^{ 1} \langle 
k^{\frac{1}{2}}_{t^{*},s^{*}, X(s^{*})} , \phi_{k} \rangle ds^{*} 
\right)^{2}dt^{*}\right] \right)\\ 
 & = \frac{1}{n}\int_0^{ 1}\int_0^{ 1} \E \left[ \left( \int_0^{ 1} \langle 
k^{\frac{1}{2}}_{t,s, X(s)} , \phi_{j} \rangle ds \right)^{2}\left(\int_0^{ 1} 
\langle k^{\frac{1}{2}}_{t^{*},s^{*}, X(s^{*})} , \phi_{k} \rangle ds^{*} 
\right)^{2}\right] dt dt^{*}.    
\end{align*}
 Using Cauchy-Schwartz inequality again
\begin{align*}
 &\E\langle \phi_{j}, (C_{n}-C) \phi_{k} \rangle^{2}_{L^{2}}\\ 
 & \leq \frac{1}{n} \int_0^{ 1}\int_0^{ 1} \E^{\frac{1}{2}} \left( \int_0^{ 
1} \langle k^{\frac{1}{2}}_{t,s, X(s)} , \phi_{j} \rangle ds \right)^{4} 
\E^{\frac{1}{2}} \left( \int_0^{ 1} \langle k^{\frac{1}{2}}_{t^{*},s^{*}, 
X(s^{*})} , \phi_{k} \rangle ds^{*} \right)^{4} dt dt^{*}.
\end{align*}
Note that we can move to the $\mbK$ inner product to obtain:
\[
 \langle k^{\frac{1}{2}}_{t,s, 
X(s)} , \phi_{k} \rangle_{L^2}
=  \langle k_{t,s, 
X(s)} , \mcL^{1/2}_k \phi_{k} \rangle_\mbK
= (\mcL^{1/2}_k \phi_k)(t,s,X(s))
\]
and $\mcL^{1/2} \phi_k$ is a function in $\mbK$, thus we can apply 
Assumption \ref{a:main}.4 to obtain
\begin{align*}
 &\E\langle \phi_{j}, (C_{n}-C) \phi_{k} \rangle^{2}_{L^{2}}
 \leq \frac{c}{n} \E \left[ \int_0^{ 1} \left ( \int_0^{ 1} \langle 
k^{\frac{1}{2}}_{t,s, X(s)} , \phi_{j} \rangle ds \right)^{2} dt\right] \E 
\left[ \int_0^{ 1} \left ( \int_0^{ 1} \langle 
k^{\frac{1}{2}}_{t^{*},s^{*}, X(s^{*})} , \phi_{k} \rangle ds^{*} \right)^{2} 
dt^{*}\right].
\end{align*}
 It is easy to see that
\begin{align*}
 &\E \left[ \int_0^{ 1} \left ( \int_0^{ 1} \langle k^{\frac{1}{2}}_{t,s, 
X(s)} , \phi_{j} \rangle ds \right)^{2} dt\right]  = \rho_{j}. 
\end{align*}
Now we obtain
\begin{equation*}
 \begin{split} 
 \E\langle \phi_{j}, (C_{n}-C) \phi_{k} \rangle^{2}_{L^{2}} & \leq c n^{-1} 
\rho_{j} \rho_{k}.
\end{split}
\end{equation*}
Therefore,
\begin{equation*}
 \begin{split} 
 \E \left(\sum_{k=1}^{\infty} \sum_{j=1}^{\infty}\frac{\rho_{k}^{-2\nu} 
\rho_{j}^{2\nu}}{(\lambda+\rho_{j})^{2}} \langle \phi_{j},(C_{n}-C) \phi_{k} 
\rangle_{L^2}\right)^{\frac{1}{2}} &\leq \left( \frac{c}{n}\sum_{k=1}^{\infty} 
\sum_{j=1}^{\infty}\frac{\rho_{k}^{1-2\nu} 
\rho_{j}^{1+2\nu}}{(\lambda+\rho_{j})^{2}} \right)^{\frac{1}{2}}.
\end{split} 
\end{equation*}
Note that 
\begin{equation*}
 \begin{split} 
 \sum_{k=1}^{\infty} \sum_{j=1}^{\infty}\frac{\rho_{k}^{1-2\nu} 
\rho_{j}^{1+2\nu}}{(\lambda+\rho_{j})^{2}} &= \sum_{k=1}^{\infty} 
\rho_{k}^{1-2\nu} \sum_{j=1}^{\infty} 
\frac{\rho_{j}^{1+2\nu}}{(\lambda+\rho_{j})^2}.
\end{split} 
\end{equation*}
Since $ 2r(1-2\nu)>1$ and $ \rho_{k}< c_{2}k^{-2r}$ 
we have
\begin{equation*}
 \begin{split} 
 \sum_{k=1}^{\infty} \rho_{k}^{1-2\nu} \leq c_{2}\sum_{k=1}^{\infty} 
k^{-2r(1-2\nu)}=c^{**}< \infty. 
\end{split} 
\end{equation*}
Finally, by applying Lemma \ref{lemma5} with $\nu_2 = 1/2$ we obtain 
\begin{equation*}
 \begin{split} 
 \E \| C^{\nu}(C+\lambda I)^{-1} (C_{n}-C) C^{-\nu} \|_{op} & \leq \gamma (n 
\lambda^{1-2\nu+\frac{1}{2r}})^{-\frac{1}{2}},
 \end{split}     
\end{equation*}
An application of Markov's inequality completes the proof.

\subsection*{Proof of Lemma \ref{lemma3}}
Recall that 
\begin{equation*}
 \begin{split} 
 \| C^{\nu_2}(C+\lambda I)^{-1} (C_{n}-C) C^{-\nu} \|^{2}_{op} & = 
\sup_{h: \| h \|_{L^2}=1} \| C^{\nu_2}(C+\lambda I)^{-1} (C_{n}-C) 
C^{-\nu} h \|^{2}_{L^2}.
\end{split} 
\end{equation*}
Note that 
\begin{equation*}
 \begin{split} 
 C^{\nu_2}(C+\lambda I)^{-1} & = 
\sum_{j=1}^{\infty}\frac{\rho_{j}^{\nu_2}}{\rho_{j}+\lambda}( \phi_{j} 
\otimes \phi_{j}),\\ 
\end{split} 
\end{equation*}
and recall that from the proof of Lemma \ref{lemma2},
\begin{equation*}
 \begin{split} 
 (C_{n}-C) C^{-\nu} h & = \sum_{k=1}^{\infty}a_{k}\rho_{k}^{-\nu} 
(C_{n}-C)\phi_{k}.
\end{split} 
\end{equation*}
Therefore
\begin{equation*}
 \begin{split} 
 C^{\nu_2}(C+\lambda I)^{-1} (C_{n}-C) C^{-\nu}h& = 
\sum_{j=1}^{\infty}\frac{\rho_{j}^{\nu_2}}{\rho_{j}+\lambda}\langle 
\sum_{k=1}^{\infty}a_{k}\rho_{k}^{-\nu} (C_{n}-C)\phi_{k},\phi_{j} 
\rangle_{L^2}\phi_{j}.
\end{split} 
\end{equation*}
By Parseval's identity we obtain 
\begin{equation*}
 \begin{split} \| C^{\nu_2}(C+\lambda I)^{-1} (C_{n}-C) 
C^{-\nu}\|^{2}_{op}& = 
\sum_{j=1}^{\infty}\left[\frac{\rho_{j}^{\nu_2}}{\rho_{j}+\lambda}\langle 
\sum_{k=1}^{\infty}a_{k}\rho_{k}^{-\nu} (C_{n}-C)\phi_{k},\phi_{j} 
\rangle_{L^2}\right]^{2}\\   
\end{split} 
\end{equation*}
By Cauchy-Schwartz inequality and using the same steps in the proof of Lemma \ref{lemma2}, 
we have 
\begin{equation*}
 \begin{split} 
 \| C^{\nu_2}(C+\lambda I)^{-1} (C_{n}-C) C^{-\nu} \|^{2}_{op} & \leq 
\sum_{k=1}^{\infty} \sum_{j=1}^{\infty}\frac{\rho_{k}^{-2\nu} 
\rho_{j}^{2\nu_2}}{(\lambda+\rho_{j})^{2}} \langle \phi_{j},(C_{n}-C) \phi_{k} 
\rangle^{2}_{L^2}.
\end{split} 
\end{equation*}
Using the same arguments in the proof of Lemma \ref{lemma2}, we obtain
\begin{equation*}
 \begin{split} 
 \E \left(\sum_{k=1}^{\infty} \sum_{j=1}^{\infty}\frac{\rho_{k}^{-2\nu} 
\rho_{j}^{2\nu_2}}{(\lambda+\rho_{j})^{2}} \langle \phi_{j},(C_{n}-C) \phi_{k} 
\rangle_{L^2}\right)^{\frac{1}{2}} 
&\leq \left( \frac{c}{n}\sum_{k=1}^{\infty} 
\sum_{j=1}^{\infty}\frac{\rho_{k}^{1-2\nu} \rho_{j}^{1+2\nu_2}}{(\lambda+\rho_{j})^{2}} 
\right)^{\frac{1}{2}}.
\end{split} 
\end{equation*}
Note that 
\begin{equation*}
 \begin{split} 
 \sum_{k=1}^{\infty} \sum_{j=1}^{\infty}\frac{\rho_{k}^{1-2\nu} 
\rho_{j}^{1+2\nu_2}}{(\lambda+\rho_{j})^{2}} &= \sum_{k=1}^{\infty} \rho_{k}^{1-2\nu} 
\sum_{j=1}^{\infty} \frac{\rho_{j}^{1+2\nu_2}}{(\lambda+\rho_{j})^2}.
\end{split} 
\end{equation*}
Note that the condition $2r(1-2\nu)>1$ implies $ \nu < 
\frac{1}{2}-\frac{1}{2r} <\frac{1}{2}$. We therefore have 
\begin{equation*}
 \begin{split} 
\left(\frac{\lambda+\rho_{j}}{\rho_{j}}\right)^{2\nu-1}& \leq 1.
\end{split} 
\end{equation*}
It follows that
\begin{equation*}
 \begin{split} 
 \sum_{k=1}^{\infty} \sum_{j=1}^{\infty}\frac{\rho_{k}^{1-2\nu} 
\rho_{j}^{1+2\nu_2}}{(\lambda+\rho_{j})^{2}}&\leq \sum_{k=1}^{\infty} 
\rho_{k}^{1-2\nu} \sum_{j=1}^{\infty} 
\frac{\rho_{j}^{2\nu+2\nu_2}}{(\lambda+\rho_{j})^{2\nu+1}}.
\end{split} 
\end{equation*}
Recall that
\begin{equation*}
 \begin{split} 
 \sum_{k=1}^{\infty} \rho_{k}^{1-2\nu} \leq c_{2}\sum_{k=1}^{\infty} 
k^{-2r(1-2\nu)}=c^{**}< \infty. 
\end{split} 
\end{equation*}
Therefore
\begin{equation*}
 \begin{split} 
 \E \| C^{\nu_2}(C+\lambda I)^{-1} (C_{n}-C) C^{-\nu} \|_{op} 
    & \leq \left( \frac{c c^{**}}{n} 
\sum_{j=1}^{\infty}\frac{ \rho_{j}^{2\nu+2\nu_2}}{(\lambda+\rho_{j})^{2\nu+1}} 
\right)^{\frac{1}{2}}.
\end{split} 
\end{equation*}
By applying Lemma \ref{lemma5} we have that
\begin{equation*}
 \begin{split} 
 \E \| C^{\frac{1}{2}}(C+\lambda I)^{-1} (C_{n}-C) C^{-\nu} \|_{op} & \leq 
\beta (n \lambda^{\frac{1}{2r} + 1 - 2\nu_2})^{-\frac{1}{2}},
 \end{split}     
\end{equation*}
where $\beta$ is a constant.  An application of the Markov inequality completes the proof.

\subsection{Proof of Lemma \ref{lemma5}}
Using the same arguments as in \cite{cai2012minimax}, we get that
\begin{align*}
\sum_{j=1}^\infty \frac{s_j^{2\nu_2 + 2\nu}}{(\lambda + s_j)^{1+2\nu}}
& = \sum_{j=1}^\infty \frac{s_j^{1 + 2\nu}}{(\lambda + s_j)^{1+2\nu}} s_j^{2\nu_2 -1} \\
& \leq \sum_{j=1}^\infty \frac{c_1^{1+2\nu}k^{-2r(1 + 2\nu)}}{(\lambda + c_2k^{-2r})^{1+2\nu}} k^{-2r(2 \nu_2 -1)} \\
& =  c_1^{1+2\nu} \sum_{j=1}^\infty \frac{k^{-2r(2 \nu_2 -1)}}{(\lambda k^{2r} + c_2)^{1+2\nu}}  \\
& \leq c_1^{1+2\nu} \left(\frac{1}{c_2} + \int_1^\infty  \frac{x^{-2r(2 \nu_2 -1)}}{(\lambda x^{2r} + c_2)^{1+2\nu}} \ dx \right)\\
& =  c_1^{1+2\nu} \left(\frac{1}{c_2} + \lambda^{2\nu_2-1 -1/2r}\int_{\lambda^{1/2r}}^\infty  \frac{y^{-2r(2 \nu_2 -1)}}{(y^{2r} + c_2)^{1+2\nu}} \ dy \right).
\end{align*}
For the integral to be finite, it is enough if $2r(2\nu_2 + 2 \nu) \geq 1+\delta$, for some $\delta>0$, as the integrand will go to zero faster than $y^{-(1+\delta)}$.  The argument for the lower bound follows the same arguments.

\clearpage
\bibliographystyle{abbrvnat}
\bibliography{references}

\end{document}